\documentclass[a4paper,11pt,final]{amsart}
\usepackage{amsfonts,amssymb,amsthm,amsmath}
\usepackage[all]{xy}

\theoremstyle{plain}

\newtheorem*{Th*}{Theorem}

\newtheorem*{Cor*}{Corollary}

\theoremstyle{definition}

\theoremstyle{remark}

\numberwithin{equation}{section}

\makeatletter
\def\Set@Scallop[#1]#2#3{{#1}\Parens{#2}{#3}}
\newcommand\DeclareScalableOperator[2]{%
  \expandafter\def\csname#1\endcsname{\@ifnextchar[{{#2}\Set@Scallop}{{#2}\Set@Scallop[{}]}}
}
\makeatother

\DeclareScalableOperator{Ct}{\mathcal{C}} % continuous functions/maps
\DeclareScalableOperator{Cc}{\mathcal{K}} % continuous w/ cpt support
\DeclareScalableOperator{Lp}{\mathbf{L}} % L^p spaces
\DeclareScalableOperator{Hom}{\mathrm{Hom}} % homomorphisms
\DeclareScalableOperator{End}{\mathrm{End}} 
\DeclareScalableOperator{Aut}{\mathrm{Aut}} 
\DeclareScalableOperator{Uenv}{\mathfrak{U}} % universal enveloping algebra
\DeclareScalableOperator{Presh}{\mathrm{Presh}} % presheaves
\DeclareScalableOperator{Sh}{\mathrm{Sh}} % sheaves
\DeclareScalableOperator{Sym}{\mathrm{Sym}} % symmetric multilinear maps
\DeclareScalableOperator{Alt}{\mathrm{Alt}} % alternating multilinear maps

\newcommand{\fa}{for all }
\newcommand{\Fs}{For some }
\newcommand{\fs}{for some }
\newcommand\mathfa[1][{}]{\quad\text{\fa{#1} }}

\newcommand{\scth}{such that }
\newcommand{\AND}{and}

\newcommand\mathtxt[1]{\quad\text{{#1}}\quad}

\newcommand{\nd}{\mathtxt\AND}

\newcommand\vphi{\varphi}
\newcommand\vrho{\varrho}

\newcommand\eps{\varepsilon}
\newcommand\nats{\mathbb{N}}
\newcommand\ints{\mathbb{Z}}
\newcommand\rats{\mathbb{Q}}
\newcommand\reals{\mathbb{R}}
\newcommand\cplxs{\mathbb{C}}
\newcommand\knums{\mathbb K}

\newcommand\vvoid{\varnothing}
\newcommand\sle{\leqslant}
\newcommand\sge{\geqslant}

\DeclareMathOperator\dom{\mathrm{dom}}

\DeclareMathOperator\id{\mathrm{id}}
\DeclareMathOperator\Top{\mathrm{Top}}
\DeclareMathOperator\Sets{\mathrm{Sets}}
\DeclareMathOperator\TopSVec{\mathrm{TopSVec}}
\DeclareMathOperator\TopSMod{\mathrm{TopSMod}}
\DeclareMathOperator\SDom{\mathrm{SDom}}
\DeclareMathOperator\SMan{\mathrm{SMan}}
\DeclareMathOperator\SRSp{\mathrm{SRSp}}

\makeatletter
\newcommand\Size[7][1]{% produces #3#4#5#6#7, where #3,#5,#7 are operators of size #2(0-4); if #1 is 0, #5 is not resized
                                 \ifx#20%
                                        \def\r@l{}\def\r@m{}\def\r@r{}%
                                 \else%
                                    \ifx#21%
                                           \def\r@l{\bigl}\def\r@r{\bigr}\def\r@m{\bigm}%
                                    \else%
                                           \ifx#22%
                                                 \def\r@l{\Bigl}\def\r@r{\Bigr}\def\r@m{\Bigm}%
                                            \else%
                                                 \ifx#23%
                                                        \def\r@l{\biggl}\def\r@r{\biggr}\def\r@m{\biggm}%
                                                  \else
                                                        \ifx#24%
                                                        \def\r@l{\Biggl}\def\r@r{\Biggr}\def\r@m{\Biggm}%
                                                        \fi%
                                                  \fi%
                                            \fi%
                                      \fi%
                                 \fi%
                                 \ifx#10%
                                       \def\r@m{}%
                                 \fi%
                                 \r@l#3{#4}\r@m#5{#6}\r@r#7%
}%

\makeatother

\newcommand\Set[3]{% set theoretical brackets #1:size (0-4), #2:elements, #3:defining predicate
                                 \Size{#1}{\{}{#2}{|}{#3}{\}}%
}%
\newcommand\Parens[2]{% Parentheses of size #1 (0-4)
  \Size[0]{#1}{(}{#2}{}{}{)}
}

\newcommand\Abs[2]{% absolute value of size #1 (0-4)
  \Size[0]{#1}{\lvert}{#2}{}{}{\rvert}
}
\newcommand\Span[2]{% span of size #1 (0-4)
  \Size[0]{#1}{\langle}{#2}{}{}{\rangle}
}
\makeatletter

\newcommand{\IfUpperCase}[1]{\begingroup 
  \protected@edef\@tempa{\expandafter\@firstofone\@firstofone#1.}%
  \expandafter\IfUpperCasE \@tempa\delimiter}

\def\IfUpperCasE #1#2\delimiter{%
  \protected@edef\@tempa{\meaning#1\meaning a}%
  \ifnum \expandafter\IfUppercaSE\@tempa \IfUppercaSE
   \endgroup \expandafter\@firstoftwo
  \else
   \endgroup \expandafter\@secondoftwo
  \fi}

\@ifundefined{strip@prefix}{\def\strip@prefix#1>{}}{}
\def\@tempa{the letter }
\edef\@tempa{\expandafter\strip@prefix\meaning\@tempa}

\expandafter\def\expandafter\IfUppercaSE\expandafter#\expandafter1\@tempa#2#3\IfUppercaSE{\uccode`#2=`#2 }

\newif\ifuc@se
\def\setuc@se#1{\IfUpperCase{#1}{\uc@setrue}{\uc@sefalse}}

\def\theoremn@me#1{\ifuc@se \lowercase{\csname#1name\endcsname}\ignorespaces%
  \else \edef\@temp{\lowercase{\lowercase{\csname#1name\endcsname}}}\@temp\ignorespaces%
  \fi}
\def\theoremn@mes#1{\ifuc@se \lowercase{\csname#1names\endcsname}\ignorespaces%
  \else \edef\@temp{\lowercase{\lowercase{\csname#1names\endcsname}}}\@temp\ignorespaces%
  \fi}
  
\def\thmref#1#2{\setuc@se{#1}\lowercase{{\theoremn@me{#1}\lowercase{\ref{#1:#2}}}}}

\newcommand{\DefTheorem}[2]{\newenvironment{#1}[1][\empty]{\ignorespaces\begin{#2}\ifx##1\empty{}\else\lowercase{\label{#1:##1}}\fi\ignorespaces}{\end{#2}\ignorespacesafterend}}

\DefTheorem{Th}{theorem}
\DefTheorem{Prop}{proposition}
\DefTheorem{Cor}{corollary}
\DefTheorem{Lem}{lemma}
\DefTheorem{Def}{definition}
\DefTheorem{Rem}{remark}
\DefTheorem{Par}{para}

\newenvironment{Defn}[1][\empty]{\ignorespaces\noindent\ignorespaces\ifx#1\empty{}\begin{Def}\else\begin{Def}[#1]\fi\ignorespaces}{\hfill$\triangle$\end{Def}\ignorespacesafterend}

\newenvironment{Par*}{\ignorespaces\noindent\ignorespaces}{\ignorespacesafterend}
\makeatletter

\newif\if@smallmat
\newif\if@none
\newif\if@paren
\newif\if@brack
\newif\if@brace
\newif\if@vline
                                 {\ifx#20%
                                        \@smallmattrue%
                                  \else%
                                         \@smallmatfalse
                                  \fi%
                                  \ifx#11%
                                         \@nonefalse\@parentrue\@brackfalse\@bracefalse\@vlinefalse%
                                  \else%
                                       \ifx#12%
                                            \@nonefalse\@parenfalse\@bracktrue\@bracefalse\@vlinefalse%
                                        \else%
                                            \ifx#13%
                                                 \@nonefalse\@parenfalse\@brackfalse\@bracetrue\@vlinefalse%
                                            \else%
                                                 \ifx#14%
                                                       \@nonefalse\@parenfalse\@brackfalse\@bracefalse\@vlinetrue
                                                 \else%
                                                       \ifx#15%
                                                             \@nonefalse\@parenfalse\@brackfalse\@bracefalse\@vlinefalse%
                                                       \else%
                                                             \@nonetrue\@parenfalse\@brackfalse\@bracefalse\@vlinefalse%
                                                       \fi%
                                                 \fi%
                                            \fi%
                                        \fi%
                                   \fi%
                                   \if@smallmat%
                                        \if@none%
                                             \begin{smallmatrix}%
                                        \else%
                                            \if@paren%
                                                  \bigl(\begin{smallmatrix}%
                                            \else%
                                                  \if@brack%
                                                          \bigl[\begin{smallmatrix}%
                                                  \else%
                                                          \if@brace%
                                                               \bigl\{\begin{smallmatrix}%
                                                          \else%
                                                               \if@vline%
                                                                    \bigl\lvert\begin{smallmatrix}%
                                                                \else%
                                                                    \bigl\lEert\begin{smallmatrix}%
                                                                \fi%
                                                          \fi%
                                                  \fi%
                                            \fi%
                                        \fi%
                                   \else%
                                        \if@none%
                                             \begin{matrix}%
                                        \else%
                                            \if@paren%
                                                  \begin{pmatrix}%
                                            \else%
                                                  \if@brack%
                                                          \begin{bmatrix}%
                                                  \else%
                                                          \if@brace%
                                                               \begin{Bmatrix}%
                                                          \else%
                                                               \if@vline%
                                                                    \begin{vmatrix}%
                                                                \else%
                                                                    \begin{Ematrix}%
                                                                \fi%
                                                          \fi%
                                                  \fi%
                                            \fi%
                                        \fi%
                                   \fi}%
                                  {\if@smallmat%
                                        \if@none%
                                             \end{smallmatrix}%
                                        \else%
                                            \if@paren%
                                                  \end{smallmatrix}\bigr)%
                                            \else%
                                                  \if@brack%
                                                          \end{smallmatrix}\bigr]%
                                                  \else%
                                                          \if@brace%
                                                               \end{smallmatrix}\bigr\}%
                                                          \else%
                                                               \if@vline%
                                                                    \end{smallmatrix}\bigr\rvert%
                                                                \else%
                                                                    \end{smallmatrix}\bigr\rEert%
                                                                \fi%
                                                          \fi%
                                                  \fi%
                                            \fi%
                                         \fi%
                                   \else%
                                        \if@none%
                                             \end{matrix}%
                                        \else%
                                            \if@paren%
                                                  \end{pmatrix}%
                                            \else%
                                                  \if@brack%
                                                          \end{bmatrix}%
                                                  \else%
                                                          \if@brace%
                                                               \end{Bmatrix}%
                                                          \else%
                                                               \if@vline%
                                                                    \end{vmatrix}%
                                                                \else%
                                                                    \end{Ematrix}%
                                                                \fi%
                                                          \fi%
                                                  \fi%
                                            \fi%
                                        \fi%
                                   \fi}%

\makeatother

\def\ger{\mathfrak}

\title[Infinite-dimensional supermanifolds]{Infinite-dimensional supermanifolds\\ over arbitrary base fields}

\author[A.~Alldridge]{Alexander Alldridge}
\address{Institut f\"ur Mathematik\\ Universit\"at Paderborn\\ Warburger Stra\ss{}e~100\\ 33100 Paderborn}
\email{alldridg@math.upb.de}

\author[M.~Laubinger]{Martin Laubinger}
\address{Institut f\"ur Mathematik\\ Universit\"at Paderborn\\ Warburger Stra\ss{}e~100\\ 33100 Paderborn}
\email{mlaubing@math.upb.de}



\thanks{This research was supported partly by the IRTG ``Geometry and Analysis of Symmetries'' grant, funded by Deutsche Forschungsgemeinschaft (DFG), Minist\`ere de l'\'Education Nationale (MENESR), and Deutsch-Franz\"osische Hochschule (DFH-UFA); and by the SFB/TR 12 ``Symmetries and Universality in Mesoscopic Systems'' grant, funded by Deutsche Forschungsgemeinschaft (DFG)}
\begin{document}

\maketitle

\begin{abstract}
	In his recent investigation of a super Teichm\"uller space, Sachse \cite{sachse-diss}, based on work of Molotkov \cite{molotkov}, has proposed a theory of Banach supermanifolds using the `functor of points' approach of Bernstein and Schwarz. We prove that the the category of Berezin--Kostant--Leites supermanifolds is equivalent to the category of finite-dimensional Molotkov--Sachse supermanifolds. Simultaneously, using the differential calculus of Bertram--Gl\"ockner--Neeb \cite{BGlN}, we extend Molotkov--Sachse's approach to supermanifolds modeled on Hausdorff topological super-vector spaces over an arbitrary non-dis\-crete Hausdorff topological base field of characteristic zero. We also extend to locally $k_\omega$ base fields the `DeWitt' super\-mani\-folds considered by Tuynman in his monograph \cite{tuynman}, and prove that this leads to a category which is isomorphic to the full subcategory of Molokov--Sachse supermanifolds modeled on locally $k_\omega$ spaces. 
\end{abstract}

\section{Introduction}

Supermanifolds were first introduced in the 1970s in an effort to find a conceptual framework in which commuting and anticommuting variables could be treated on an equal footing. The physical motivation stems from quantum field theory in its functional integral formulation, which describes fermionic particles by anticommuting fields. Since their introduction, there have been several approaches to the rigorous mathematical definition of supermanifolds (and their morphisms). The most commonly used (in terms of ringed spaces) is due to Berezin and Leites \cite{berezin, leites}. Kostant \cite{kostant-supergeom} proposed an approach via Hopf algebras which is equivalent if the base field is $\reals$. The definition of supermanifolds which is probably closest to the physicist's usage of the concept was given by DeWitt \cite{dewitt}, and in his monograph \cite{tuynman}, Tuynman provides a rigorous account of this theory from first principles. 

In the physical applications, not only finite-dimensional, but also \emph{infinite-dimensional} supermanifolds arise naturally, \emph{e.g.}~as direct limits. Other natural occurrences of infinite-dimensional supermanifolds would be, \emph{e.g.}, mapping spaces (for instance, gauge supergroups) and supergroups of superdiffeomorphisms. 

Unfortunately, in infinite dimensions, the abundance of definitions advocated by different authors is at least as large as in finite dimensions, and there seems to be even less agreement as to which of these is to be favoured. For example, Schmitt has proposed a definition which generalises the ringed space approach \cite{schmitt-infdimsmf}, in order to define infinite-dimensional real and complex \emph{analytic} supermanifolds. This has the benefit of using a well-established conceptual framework; however, it is not possible to define, in full generality, an internal $\mathrm{Hom}$---\emph{i.e.}~there is, for supermanifolds $\mathcal M$ and $\mathcal N$, in general no supermanifold of morphisms $\mathcal M\to\mathcal N$ in Schmitt's approach. Kostant's approach has also been extended to infinite dimensions by Jask\'olski \cite{jaskolski-smoothcoalg}; it is limited to the case of the base field $\reals$. 

Both obstacles (limitation to the real field and non-existence of an internal $\mathrm{Hom}$) were overcome by Molotkov \cite{molotkov}. His idea is to define Banach supermanifolds via their functor of higher points. This approach goes back, in the finite-dimensional setting, to Bernstein \cite{leites,manin,deligne-morgan}. It is an idea due to Schwarz \cite{schwarz} to use this as the \emph{definition} of supermanifolds. Although in finite dimensions, supermanifolds are not usually defined by this method, the `point functor' is nonetheless an indispensable tool, \emph{e.g.}, to introduce and study the general linear supergroup, and the generalisations of the other classical groups, to the super framework \cite{manin,schmitt-supergeom}. 

Many of the details of Molotkov's approach were worked out by Sachse in his Ph.D.~thesis \cite{sachse-diss}. He uses this framework to good effect in the study of a super Teichm\"uller space. In joint work with Wockel \cite{sachse-wockel}, he has used it to define the superdiffeomorphism supergroup of a compact supermanifold. This appears to be evidence in favour of Molotkov--Sachse's approach. 

\smallskip\noindent
Given its level of abstraction, it appears appropriate to explain the idea behind the \emph{functor of points} approach followed by Molotkov--Sachse. 

In its simplest form, a set $X$ can be understood as a collection of points. A \emph{point} is just a map $*\to X$ where $*$ is some singleton set, fixed once and for ever. If now $f:X\to Y$ is a map of sets, then the equation $y=f(x)$ is the same as the commutative diagram
\[
	\xymatrix@C+15pt@R-5pt{%
		&{*}\ar[ld]\ar[rd]\\
		X\ar[rr]_f&&Y
	}
\]
Thus, the map $f$ can be understood as the collection of all such diagrams.

Let us now switch to supermanifolds, defined as usual as graded ringed spaces which are locally isomorphic to $(\reals^p,\mathcal C^\infty_{\reals^p}\otimes\bigwedge(\reals^q)^*)$. Assume that $Z$ is a supermanifold whose underlying space is $*$. Then $Z=(*,\lambda)$ where $\lambda$ is some (finitely generated) Grassmann algebra; if $Y=(Y_0,\mathcal O_Y)$, then any morphism $h:Z\to Y$ is determined by a pair $(h(*),h^*:\mathcal O_{Y,h(*)}\to\lambda)$ where $h(*)\in Y_0$ and $h^*:\mathcal O_Y(Y_0)\to\lambda$ is an even algebra morphism on the global sections module of $\mathcal O_Y$. There may be a variety of such morphisms, so $h$ is \emph{not} characterised solely by its value $h(*)$.

So, in place of $*$, we consider all of these the ringed spaces $Z$ as above, \emph{i.e.}~the spaces $*_\lambda=(*,\lambda)$ where $\lambda$ is \emph{any} finitely generated Grassmann algebra. (These are exactly the \emph{higher points} we referred to earlier.) A morphism $f:X\to Y$ is uniquely determined by the commutative diagrams
\[
	\xymatrix@C+15pt@R-5pt{%
		&{*_\lambda}\ar[ld]\ar[rd]\\
		X\ar[rr]_{f_\lambda}&&Y
	}
\]
More precisely, the collections $(*_\lambda\to X)_\lambda$ and $(*_\lambda\to Y)_\lambda$ are functors defined on the category of Grassmann algebras, which take values in the category of sets. Moreover, the family $(f_\lambda)$ determining $f$ is comprised of maps $\mathrm{Mor}(*_\lambda,X)\to\mathrm{Mor}(*_\lambda,Y)$, and defines a `functor morphism' or `natural transformation', since it behaves naturally under morphisms in the category of Grassmann algebras. 

Clearly, not every natural transformation comes from a morphism. To understand when it does, assume that $X=\reals^{m|n}$ and $Y=\reals^{p|q}$. We observe that the set $\underline X(\lambda)$ of morphisms $*_\lambda\to X$ is exactly $(\reals^{n|m}\otimes\lambda)_0$, the even part of the tensor product $\reals^{m|n}\otimes\lambda$. (Here, $\reals^{m|n}$ denotes the graded vector space $\reals^m\oplus\reals^n$.) In particular, it has the structure of a finite-dimensional vector space over $\reals$, and of a module over the commutative $\reals$-algebra $\lambda_0$. 

Then $(f_\lambda)$ is defined by a (unique) morphism $f$ if and only if every map
\begin{equation}\label{eq:ms-mordef}
	f_\lambda:(\reals^{m|n}\otimes\lambda)_0\to(\reals^{p|q}\otimes\lambda)_0
\end{equation}
is smooth, and its first derivative $df_\lambda(x)$ at any $x$ is $\lambda_0$-linear. Molotkov--Sachse's idea (due to Schwarz in the finite-dimensional case) is now to use this as the \emph{definition} of morphisms of superdomains, replacing $\reals^{m|n}$, $\reals^{p|q}$ by any pair of graded Banach spaces. At this point, it is necessary to consider the above functor as taking values in the category of \emph{topological spaces} (and not only of sets). Supermanifolds and their morphisms are then built up by using an appropriate generalisation of the concept of `local charts' in the categorical framework. (Technically, one uses Grothendieck topologies.)

\smallskip\noindent
While it seems to be general enough for serious applications, the equivalence of this approach with others, even in finite dimensions, seems not to have been completely worked out, at least in published form. (For the base field $R=\reals$ and finite dimensions, it has been shown by Voronov \cite{voronov-maps} that the categories of local models (\emph{i.e.}~superdomains) in the Berezin--Kostant--Leites and Bernstein--Schwarz sense are equivalent, but beyond this, there does not appear to be any published reference.) 

Moreover, the limitation to Banach spaces as model spaces seems artificial. In fact, there is even no \emph{a priori} reason why other base fields (of zero characteristic) should be excluded from the mathematical study of supermanifolds. \emph{E.g.},~we mention the field of formal Laurent series, as well as the $p$-adic and adelic fields. Even in the finite-dimensional case, such a setting seems not to have been considered before. 

Bertram--Gl\"ockner--Neeb \cite{BGlN} have proposed a natural and robust $\mathcal C^k$ calculus which is valid under very general assumptions on the base field---and, with some restrictions, even for base rings. This approach generalises the usual differential calculi for ultrametric fields, and for the real and complex fields; it does away with assumptions such as completeness and local convexity of model spaces in the real and complex cases. 

Besides the implications this has for the generality and scope of a theory based on these ideas, it allows without ado for a differential calculus over Grassmann algebras. This leads to a rigorous formulation of the intuitive idea that a morphism of superdomains is a map which is differentiable over a Grassmann algebra, and therefore conceptually simplifies some aspects of the Molotkov--Sachse calculus. 

\smallskip\noindent
To explain our generalisation of Molotkov--Sachse's idea, let us recall the $\mathcal C^k$ calculus of Bertram--Gl\"ockner--Neeb. Consider a function $f:U\to\reals^n$ defined on an open subset $U\subset\reals^m$. Then $f$ is $\mathcal C^1$ (over $\reals$) if and only if there exists a continuous map $f^{[1]}$ \scth
\begin{equation}\label{eq:hadc1}
	f(x+tv)-f(x)=t\cdot f^{[1]}(x,v,t)
\end{equation}
where $f^{[1]}$ is defined on an appropriate open subset $U^{[1]}\subset U\times\reals^m\times\reals$. 

The key realisation of Bertram--Gl\"ockner--Neeb is that this approach, which insists on separating the variables of all directional derivatives, gives a useful definition of $\mathcal C^1$ maps (and, by induction, of $\mathcal C^k$ maps for $k\in\nats\cup\{\infty\}$) if one replaces $\reals$ by a suitable topological ring $R$, and $\reals^m$, $\reals^n$ by Hausdorff topological $R$-modules.\footnote{Here, `suitable' means that $R$ is Hausdorff and has a dense group of units; for all but the most basic questions of differential calculus, one also has to assume that $R$ be commutative.} As a matter of terminology, we say a map $f$ satisfying a condition analogous to \eqref{eq:hadc1} is $\mathcal C^1$ \emph{over $R$}. 

If $R$ is such a ring, and $\lambda=R[\theta_1,\dotsc,\theta_n]$ is any finitely generated Grassmann algebra over $R$, then the even part $\lambda_0$ also satisfies the conditions necessary for Bertram--Gl\"ockner--Neeb's definition of $\mathcal C^k$ to be applicable over the ring $\lambda_0$ (for a suitable topology on $\lambda$). Hence, it seems natural to use the $\lambda_0$-module structure on $(E\otimes\lambda)_0$, and to define a natural transformation given by a collection $(f_\lambda)$, $f_\lambda:(E\otimes\lambda)_0\to(F\otimes\lambda)_0$ to be $\mathcal C^k_{MS}$ if each $f_\lambda$ is $\mathcal C^k$ over~$\lambda_0$. 

More precisely, fix a base ring $R$. If $E$, $F$ are graded Hausdorff topological $R$-modules, let $\underline E$ be the functor with values in topological spaces given by $\underline E(\lambda)=(E\otimes\lambda)_0$ (with a suitable topology). If $\mathcal U$ is a functor \scth $\mathcal U(\lambda)\subset\underline E(\lambda)$ is open \fa $\lambda$, and $f:\mathcal U\Rightarrow\underline F$ is a natural transformation, we say that $f$ is $\mathcal C^1_{MS}$ if there exists a natural transformation $f^{[1]}:\mathcal U^{[1]}\Rightarrow\underline F$ (where $\mathcal U^{[1]}\subset\underline{E\times E\times R}$ is suitably defined), \scth 
\[
	f_\lambda(x+tv)-f_\lambda(x)=t\cdot f_\lambda^{[1]}(x,v,t)\mathfa\lambda\in\Lambda\,,\,(x,v,t)\in\mathcal U^{[1]}(\lambda)\ .
\]
(This can also be phrased in terms of categorified linear algebra, see \thmref{Def}{c1mordef}.) Although this is not obvious, it turns out that $\mathcal C^\infty$ natural transformations can then indeed be equivalently characterised by Molotkov--Sachse's smoothness condition~\eqref{eq:ms-mordef} (\emph{cf}.~\thmref{Th}{smoothequiv}). The key observation here is that $\mathcal C^\infty_{MS}$ morphisms admit an exact Taylor expansion, \emph{i.e.}, they are `Grassmann analytic' in the sense of Berezin \cite{berezin}. 

\medskip\noindent
Using this approach, we define a category $\SDom_{MS}$ of Molotkov--Sachse superdomains which serve as the local models for a category $\SMan_{MS}$ of (possibly infinite-dimensional) Molotkov--Sachse supermanifolds over any non-discrete Hausdorff topological field $R$ of characteristic zero. We also define a category $\SMan_{BKL}$ of (finite-dimensional) Berezin--Kostant--Leites supermanifolds over $R$, and a category of (possibly infinite-dimensional) DeWitt--Tuynman supermanifolds over $R$. 

Whereas Molotkov--Sachse supermanifolds can be defined without restrictions on the topology of the model spaces (apart from being Hausdorff), DeWitt--Tuynman supermanifolds do not extend far beyond finite dimensions: We show that the definition can be made for model spaces which are \emph{locally $k_\omega$}; essentially, this reduces to the case of direct limits of finite-dimensional spaces. The reason for the restriction to locally $k_\omega$ spaces is that countable inductive limits of topological groups will almost never be topological groups, unless they are locally $k_\omega$ (for more details, see the main text). 

This limitation seems to indicate that while DeWitt's approach is intuitive, it is inherently restricted to rather particular classes of infinite-dimensional supermanifolds (which, nonetheless, are of interest in applications). Our main results (Theorems \ref{th:smoothequiv}, \ref{th:msbklfdequiv}, \ref{th:dwtmsiso}) are as follows. 

\begin{Th*}
	If $R=\reals$ or $\cplxs$, the full subcategory of $\SMan_{MS}$ formed by the Hausdorff Molotkov--Sachse supermanifolds locally modeled on graded Banach spaces is the category of supermanifolds as defined by Molotkov--Sachse. 
\end{Th*}

\begin{Th*}
	The category $\SMan_{BKL}$ is equivalent to the full subcategory of $\SMan_{MS}$ formed by the finite-dimensional Molotkov--Sachse supermanifolds.
\end{Th*}

\begin{Th*}
	If $R$ is locally $k_\omega$ (for instance, locally compact), then the category $\SMan_{dWT}$ is isomorphic to the full subcategory of $\SMan_{MS}$ formed by the supermanifolds modeled on locally $k_\omega$ spaces. 
\end{Th*}

The corresponding statements in the first theorem for Berezin--Kostant--Leites and (finite-dimensional) DeWitt--Tuynman supermanifolds (in case $R\in\{\reals,\cplxs\}$ and $R=\reals$, respectively), \emph{i.e.}~that they generalise the definitions extant in the literature, are also correct; the proof is easy once the first theorem has been established. 

The proof of the second theorem accounts for about half of the volume of our paper, so it is perhaps appropriate to comment briefly upon the underlying idea. The proof consists of two steps: The first is to define an equivalence $\Phi:\SDom_{MS}^{fd}\to\SDom_{BKL}$ of the category of finite-dimensional Molotkov--Sachse superdomains and the category of Berezin--Kostant--Leites supermanifolds. This is done by `standard procedures'. \emph{E.g.}, the proof of faithfulness amounts to the generalisation of the identification $\mathrm{Mor}(*_\lambda,\reals^{m|n})\cong(\reals^{m|n}\otimes\lambda)_0$ to arbitrary base fields (\thmref{Prop}{lambdapoints}). 

The second step is to globalise this correspondence by `duality'. Indeed, one may consider the categories $\Sh0{\mathcal C_1}$ and $\Sh0{\mathcal C_2}$ of sheaves on $\mathcal C_1=\SDom_{MS}^{fd}$ and $\mathcal C_2=\SDom_{BKL}$, respectively. These are certain subcate\-gories of the functor categories $\Sets^{\mathcal C_j^{op}}$ , defined by specifying Grothendieck topologies on the categories $\mathcal C_j$. There is a `transpose' $\Phi^*:\Sets^{\mathcal C_2^{op}}\to\Sets^{\mathcal C_1^{op}}$ given by composition with the functor $\Phi$, and this is again an equivalence of categories. In fact, it restricts to an equivalence $\Sh0{\mathcal C_2}\to\Sh0{\mathcal C_1}$. Now, $\SMan_{MS}^{fd}$ can be embedded into $\Sh0{\mathcal C_1}$ (by the Yoneda embedding), and similarly for $\SMan_{BKL}$. Finally, one identifies the essential images of these embeddings and sees that $\Phi^*$ restricts to an equivalence between~them. 

\medskip\noindent
Let us sketch the contents of the paper. In Section \ref{section2}, we define the functorial framework for Molotkov--Sachse superdomains over an arbitrary base field and discuss the notion of smoothness for natural transformations. The results of this section are basic for all that follows; in particular, we show in \thmref{Th}{smoothequiv} that our definition of smoothness is the same as Molotkov--Sachse's. 

In Section \ref{section31}, we define Berezin--Kostant--Leites superdomains over an arbitrary base field, and prove the equivalence of this category with the category of finite-dimensional Molotkov--Sachse superdomains in \thmref{Th}{msbklsdomequiv}. The constructions in this section are parallel to those in the case of the base field $\reals$. In Sections \ref{section32} and \ref{section33}, we globalise the equivalence. Our approach is to embed both categories of supermanifolds into appropriate categories of sheaves on superdomains, for some correctly defined Grothendieck topologies. By dual\-ity, we get an equivalence of the categories of sheaves, and this restricts to an equivalence of the categories of supermanifolds. 

One obtains along the way a criterion (\thmref{Prop}{msrepresentable}) for sheaves on superdomains to be representable as supermanifolds, which might be of independent interest. (Another instance of a representability criterion, of a slightly different flavour, is to be found in \cite{fioresi-lledo-varadarajan}.) Let us remark that using these ideas it should be possible to rigorously define the correct notions of smooth superstacks and in particular, superorbifolds, and their morphisms. This might be useful in the study of supermoduli problems. 

Finally, in Section \ref{section4}, we show that DeWitt--Tuynman's definition of supermanifolds can be simplified, and at the same time extended to arbitrary base fields $R$ and infinite dimensions. Using our results on Molotkov--Sachse's definition of supermani\-fold morphisms from Section \ref{section2}, this again gives rise to a category of supermanifolds which is \emph{isomorphic} to the full subcategory of $\SMan_{MS}$ of supermanifolds modeled on locally $k_\omega$ spaces (\thmref{Th}{dwtmsiso}).

\medskip\noindent
\emph{Acknowledgements.} We wish to thank the anonymous referee for the diligent reading of the paper, and precise comments on its overall readability. Should it have been improved, then this is entirely to his or her credit. The first named author wishes to thank Helge Gl\"ockner for helpful remarks and references which substantially enhanced an earlier version of the manuscript. 

\section{Smooth functor morphisms}\label{section2}

We introduce the category of functors from Grassmann algebras to topological spaces. Using a categorified version of the Bertram--Gl\"ockner--Neeb differential calculus, we define smoothness for natural transformations. 

\subsection{The Grassmann category}

We begin by introducing suitable topologies on Grassmann algebras and their modules. 

\begin{Par}
	In all what follows, let $R$ be a (unital) commutative Hausdorff topological ring whose \emph{group of units} $R^\times$ is dense. We will be mainly interested in the case that $R=\knums$ is a non-discrete topological field, but to neatly construct tangent objects, one may want to consider such rings as $R=\knums[\eps]/(\eps^2)$. 
	
	We consider the supercommutative $R$-superalgebras $\lambda^n=R[\theta_1,\dotsc,\theta_n]$ and $\lambda^\infty=R[\theta_i|i\in\nats]$ freely generated by odd indeterminates $\theta_i$. On the finite-rank Grassmann algebras $\lambda^n=\prod_{\Abs0I\sle n}R\theta_I$, we consider the product topology. (Here, for $I=(i_1<\dotsc<i_k)$, we set $\theta_I=\theta_{i_1}\dotsm\theta_{i_k}$.) On $\lambda^\infty$, we consider the direct limit topology with respect to the canonical embeddings $\lambda^n\to\lambda^\infty$. This topology is Hausdorff. If $R$ is metrisable, then $\lambda^\infty$ is a topological ring if and only if $R$ is locally compact \cite[Theorems 2-4]{yamasaki}. More generally, it is a topological ring if $R$ is locally $k_\omega$ \cite[Proposition 4.7]{ggh-kac} (for instance, if $R$ is locally compact). 
	
	Here, we recall that a Hausdorff topological space $X$ is $k_\omega$ if $X=\varinjlim_{n\in\nats}K_n$ as a topological space, \fs compact subsets $K_n\subset K_{n+1}\subset X$. Moreover, a Hausdorff space $X$ is \emph{locally $k_\omega$} if every point has an open neighbourhood which is a $k_\omega$ space. Compare \cite[\S~4]{ggh-kac} for an excellent discussion of $k_\omega$ and locally $k_\omega$ spaces. Note that a locally $k_\omega$ space which is metrisable is automa\-tically locally compact \cite[Proposition 4.8]{ggh-kac}. 
	
	Let $\Lambda^\infty$ be the category whose objects are the algebras $\lambda^n$, $\lambda^\infty$, and whose morphisms are the even unital $R$-algebra morphisms $\lambda^n\to\lambda^m$. Then $\lambda^0=R$ is the null object, and the unique morphisms $\lambda\to R$ and $R\to\lambda$ will be denoted by $\eps$ and $\eta$, respectively. We let $\lambda^+=\ker\eps$ for any $\lambda\in\Lambda^\infty$. The full subcategory whose objects are $\lambda^n$ ($n<\infty$) will be denoted by $\Lambda$.  
	
	Of course, $\eta:R\to\lambda$ is continuous \fa $\lambda\in\Lambda^\infty$, and $\eps:\lambda\to R$ is continuous for $\lambda\in\Lambda$. It is also continuous for $\lambda=\lambda^\infty$, as $\eps=\varinjlim_n\eps|_{\lambda^n}$. 
	
	The following simple observation will be fundamental to our study of differential calculus over Grassmann algebras. 
\end{Par}

\begin{Lem}[lambdaunitdense]
	Let $\lambda=\lambda^n$ where $n\in\nats\cup\infty$. Then $\lambda^\times$ is dense in $\lambda$.  
\end{Lem}

\begin{proof}
	If $r\in R$, then $r-\eta\eps(r)$ is nilpotent. Hence, $\lambda^\times=\eps^{-1}(R^\times)$. Let $x\in\lambda\setminus\lambda^\times$ and $r=\eps(x)$. There exists a net $(r_\alpha)$ in $R^\times$ converging to $r$, and $x=x-\eta(r)+\lim_\alpha\eta(r_\alpha)$. Since $\eps(x-\eta(r)+\eta(r_\alpha))=r_\alpha\in R^\times$, $x\in\overline{\lambda^\times}$. 
\end{proof}

\begin{Par}
	A \emph{graded topological $R$-module} is a direct sum $E=E_0\oplus E_1$ of two topological $R$-modules. The category whose objects are the \emph{Hausdorff}Êgraded topological $R$-modules, and whose morphisms are the even continuous $R$-linear maps, will be denoted by $\TopSMod_R$. We futher denote by $\Top$ the category of topological spaces and continuous maps, and by $\Top^\Lambda$ ($\Top^{\Lambda^\infty}$) the category of functors $\Lambda\to\Top$ ($\Lambda^\infty\to\Top$) and their natural transformations. If $A$ is a (not necessarily unital) $R$-superalgebra, we write 
	\[
	\underline E(A)=(E\otimes A)_0=E_0\otimes A_0\oplus E_1\otimes A_1\ .
	\]
	
	  Let $E$ be a graded topological $R$-module. If $N$ is a non-negative integer, then $E\otimes\lambda^N=\prod_{\Abs0I\sle N}E\theta_I$. We endow $E\theta_I$ with the topology turning the canonical bijection $E\theta_I\to E$ into a homeomorphism, and take the product topology on $E\otimes\lambda^N$. We let $E\otimes\lambda^\infty=\varinjlim_NE\otimes\lambda^N$ in $\Top$. This turns $\underline E(\lambda)=(E\otimes\lambda)_0$ and $E\otimes\lambda$ into Hausdorff topological spaces, for any $\lambda\in\Lambda^\infty$. (If $x,y\in E\otimes\lambda^N$, $x\neq y$, and $U,V\subset E\otimes\lambda^N$ are open and disjoint, $x\in U$, $y\in V$, then $U^\infty=U+\sum_{\Abs0I>N}E\theta_I$, $V^\infty=V+\sum_{\Abs0I>N}E\theta_I$ are open in $E\otimes\lambda^\infty$, and $U^\infty\cap V^\infty=\vvoid$.) We call this topology the \emph{standard topology}. The standard topology is a $\lambda_0$-module topology on $\underline E(\lambda)$ and a $\lambda$-module topology on $E\otimes\lambda$ if $\lambda\in\Lambda$; for $\lambda=\lambda^\infty$, it is a $\lambda_0$-module (resp.~$\lambda$-module) topology if $R$ and $E$ are locally $k_\omega$ spaces. This follows readily from \cite[Proposition 4.7 or Corollary 5.7]{ggh-kac}. If $R$ is locally $k_\omega$ (\emph{e.g.}, locally compact), then $E$ is locally $k_\omega$ if it is the direct limit, as a topological space, of a sequence of finitely generated $R$-submodules.
	  
	Any morphism $\phi:\lambda\to\lambda'$ in $\Lambda^\infty$ gives rise to a continuous $R$-linear map $\underline E(\phi)=(\id_E\otimes\,\phi)|_{\underline E(\lambda)}:\underline E(\lambda)\to\underline E(\lambda')$. (Although for $\lambda=\lambda^\infty$, $\underline E(\lambda)$ is not in general a topological $R$-module, it is an $R$-module and topological space.) One checks that this defines functors $\underline E$ in $\Top^\Lambda$ and $\Top^{\Lambda^\infty}$. 
\end{Par}

\subsection{The DeWitt topology} We introduce a new topology on graded topological modules over Grassmann algebras. Its main purpose will be to single out certain subsets (namely, the DeWitt open subsets) which will be our generalisation of superdomains. 

\begin{Par}[dewitt-top]
	Consider the projection $\underline E(\eps):\underline E(\lambda)\to E_0$ (resp.~$\id\otimes\,\eps:E\otimes\lambda\to E$) whose kernel is $\underline E(\lambda)^+=(E\otimes\lambda^+)_0$ (resp.~$E\otimes\lambda^+$). We define the \emph{DeWitt topology} on $\underline E(\lambda)$ (resp.~$E\otimes\lambda$) to be the coarsest topology for which $\underline E(\eps)$ (resp.~$\id\otimes\,\eps$) is continuous. We will indicate its application by the subscript ${}_{dW}$. The following lemma is easy to check by hand.	
\end{Par}

\begin{Lem}[dewitt-top]
	Let $E\in\TopSMod_R$ and $\lambda\in\Lambda^\infty$. Then $U\subset\underline E(\lambda)_{dW}$ (resp.~$U\subset E\otimes\lambda$) is open if and only if $U=\underline E(\eta)(V)+\underline E(\lambda)^+$ (resp.~$U=V\otimes1+E\otimes\lambda^+$) \fs open $V\subset E_0$ (resp.~$V\subset E$). The DeWitt topology on $\underline E(\lambda)$ (resp.~$E\otimes\lambda$) is a $\lambda_0$-module (resp. $\lambda$-module) topology. \end{Lem}

\begin{proof}
	The open subsets of $\underline E(\lambda)_{dW}$ are exactly $\underline E(\eps)^{-1}(V)$, $V\subset E_0$ open. It is clear that $\underline E(\eps)^{-1}(V)=\underline E(\eta)(V)+\ker\underline E(\eps)=\underline E(\eta)(V)+\underline E(\lambda)^+$. 
	
	This collection of subsets is invariant under translations, seeing that $\underline E(\lambda)=\underline E(\eta)(E_0)\oplus\underline E(\lambda)^+$ as an $R$-module. Let $U=\underline E(\eta)(V)+\underline E(\lambda)^+$ be a DeWitt-open $0$-neighbourhood. Then $V\subset E_0$ is an open $0$-neighbourhood and there exists an open $0$-neighbourhood $V'\subset E_0$ \scth $V'+V'\subset V$. For $U'=\underline E(\eta)(V')+\underline E(\lambda)^+$, it follows that $U'+U'\subset U$, so $\underline E(\lambda)_{dW}$ is a topological monoid. 
	
	 Next, we check that the multiplication map $\lambda_0\times\underline E(\lambda)_{dW}\to\underline E(\lambda)_{dW}$ is continuous. Let $U=\underline E(\eta)(V)+\underline E(\lambda)^+$ be DeWitt-open. There exist open subsets $W'\subset R$ and $V'\subset E_0$ \scth $W'V'\subset V$. Let $W=\eps^{-1}(W')$ and $U'=\underline E(\eta)(V')+\underline E(\lambda)^+$. For any $r\in W$, we may write $r=r_0+r'$ where $r_0\in W'$ and $r'\in\lambda_0^+$. Then 
	\[
		rU'=\underline E(\eta)(r_0V')+r'\underline E(\eta)(V')+r\underline E(\lambda)^+\subset\underline E(\eta)(V)+\underline E(\lambda)^+=U\ ,
	\]
	and this proves the assertion. 
\end{proof}

\begin{Cor}
	Let $E\in\TopSMod_R$. Then $(E\otimes\lambda^\infty)_{dW}=\varinjlim_n(E\otimes\lambda^n)_{dW}$ and $\underline E(\lambda^\infty)_{dW}=\varinjlim_n\underline E(\lambda^n)_{dW}$ in $\Top$.
\end{Cor}

\begin{Defn}
	Given $F\in\Top^\Lambda$ ($F\in\Top^{\Lambda^\infty}$) and an open subset $U\subset F(R)$, define the functor $F_U\in\Top^{\Lambda^\infty}$ by 
	\[
		F_U(\lambda)=F(\eps)^{-1}(U)\ ,\ F_U(\alpha:\lambda\to\lambda')=F(\alpha)|_{F_U(\lambda)}\ .
	\]
	We call $F_U$ a \emph{restriction} of $F$. For $F=\underline E$, one obtains
	\[
		\underline E_U(\lambda)=\underline E(\eps)^{-1}(U)=U\times(E\otimes\lambda^+)_0\ ,\ \underline E_U(\alpha:\lambda\to\lambda')=\underline E(\alpha)|_{\underline E_U(\lambda)}\ .
	\]	

	Given $F,F'\in\Top^\Lambda$ ($F,F'\in\Top^{\Lambda^\infty}$), then $F'$ is called a \emph{subfunctor} of $F$ if $F'(\lambda)\subset F(\lambda)$ \fa $\lambda\in\Lambda$ and these inclusions define a natural transformation $F'\Rightarrow F$. Moreover, a subfunctor $F'\subset F$ is called \emph{open} if \fa $\lambda\in\Lambda$, $F'(\lambda)$ is open in $F(\lambda)$. 
\end{Defn}

The DeWitt topology characterises open subfunctors, as follows.

\begin{Prop}[dewittgrothtop]
	Let $E\in\TopSMod_R$, $\lambda\in\Lambda^\infty$, $U\subset\underline E(\lambda)$ be a subset. There exists an open subfunctor $\mathcal U\subset\underline E$ \scth $\mathcal U(\lambda)=U$ if and only if $U$ is DeWitt-open. Such a functor is unique, and given by $\mathcal U=\underline E_{U_R}$ where $U_R=\underline E(\eps)(U)=\mathcal U(R)$. 
\end{Prop}

\begin{proof}
	Let $\mathcal U$ be an open subfunctor of $\underline E$. By \cite[Proposition~3.5.8]{sachse-diss}, one has $\mathcal U=\underline E_V$ where $V=\mathcal U(R)$. In particular, 
	\[
		\mathcal U(\lambda)=\underline E(\eps)^{-1}(V)=E(\eta)(V)+\underline E(\lambda)^+\ .
	\]
	The assertion follows from \thmref{Lem}{dewitt-top}. 
\end{proof}

\subsection{$\mathcal C^\infty$ morphisms in $\Top^\Lambda$}

\begin{Par*}
	In this section, we introduce and study a notion of smoothness in $\Top^\Lambda$ which is a generalisation of one due to Molotkov--Sachse. We will also give a consistent definition of $\mathcal C^k$ morphisms for $k$ finite. 
	
	Recall the following general and robust notion of continuous differentiability due to Bertram--Gl\"ockner--Neeb.
\end{Par*}

\begin{Defn}[c1def]\cite{BGlN}
	Let $E,F$ be Hausdorff topological $R$-modules, $\vvoid\neq U\subset E$ be open, and $f:U\to F$ be continuous. Define $E^{[1]}=E\times E\times R$ and $U^{[1]}=\Set1{(x,v,t)\in E^{[1]}}{x,x+vt\in U}$. Then $f$ is $\mathcal C^1$ if there exists a continuous map $f^{[1]}:U^{[1]}\to F$ \scth 
	\begin{equation}\label{eq:diffdef}
		f(x+vt)-f(x)=t\cdot f^{[1]}(x,v,t)\mathfa(x,v,t)\in U^{[1]}\ .
	\end{equation}
	Inductively, define $f$ to be $\mathcal C^{k+1}$ ($k\sge1$) if $f^{[1]}$ is $\mathcal C^k$. Then $f$ is called \emph{smooth} or $\mathcal C^\infty$ if it is $\mathcal C^k$ for any $k$. 
	
	Whenever we want to stress the dependence on the base ring, we will say that $f$ is $\mathcal C^k$ (or smooth, for $k=\infty$) \emph{over $R$.}
\end{Defn}

\begin{Rem}	
	The map $f^{[1]}$ is unique, and $df(x)v=\partial_vf(x)=f^{[1]}(x,v,0)$ is $R$-linear in $v$. The higher derivatives $d^kf(x)(v_1,\dotsc,v_k)=\partial_{v_1}\dotsm\partial_{v_k}f(x)$ (defined for $f$ of class $\mathcal C^k$) are $R$-multilinear and symmetric. For all of these statements, \emph{cf}.~\cite{BGlN}.
\end{Rem}
	
The above definition lends itself to a transposition into the framework of categorical linear algebra. 

We warn the reader that we will frequently pass from the category $\Top^\Lambda$ to the category $\Top^{\Lambda^\infty}$. In fact, any morphism $f:\mathcal U\Rightarrow\underline F$ in $\Top^\Lambda$ (where $E,F\in\TopSMod_R$ and $\mathcal U\subset\underline E$ is an open subfunctor) has a unique extension to a morphism in $\Top^{\Lambda^\infty}$. However, $\Top^{\Lambda^\infty}$ is less well behaved, in particular, we recall from above that $\underline E(\lambda^\infty)$ is not a topological $\lambda_0^\infty$-module unless $E$ and $R$ are locally $k_\omega$. For this reason, constructions based on categorical linear algebra will be performed in $\Top^\Lambda$. 

We obtain the following natural notion of continuous differentiability for natural transformations. 

\begin{Defn}[c1mordef]
	Let $E\in\TopSMod_R$. The polynomial map 
	\[
	E^{[1]}=E\times E\times R\to E:(x,v,t)\mapsto x+vt
	\]
	defines a natural transformation $\alpha^1:\underline E^{[1]}\Rightarrow\underline E$ in $\Top^\Lambda$ by
	\[
		\alpha^1_\lambda:\underline E^{[1]}(\lambda)\to\underline E(\lambda):(x,v,t)\mapsto x+vt\ .
	\]
	Here, $\underline E^{[1]}(\lambda)=(E^{[1]}\otimes\lambda)_0=(E\otimes\lambda)_0\times(E\otimes\lambda)_0\times\lambda_0$. 
	
	Next, let $\mathcal U\subset\underline E$ be an open subfunctor, so that $\mathcal U=\underline E_U$ \fs unique open $U\subset E_0$. We define $\mathcal U^{[1]}=\underline E^{[1]}_{U^{[1]}}$. Then $\alpha^1=(\alpha^1_\lambda):\mathcal U^{[1]}\Rightarrow\mathcal U$. There exist natural transformations $p_1,p_2:\mathcal U^{[1]}\Rightarrow\underline E$ and $p_3:\mathcal U^{[1]}\Rightarrow\underline R$ induced by the projections of the product $\underline E^{[1]}=\underline E\times\underline E\times\underline R$. 
	
	Let $f:\mathcal U\Rightarrow\underline F$ be a natural transformation. We say that $f$ is \emph{$\mathcal C^1_{MS}$} if there exists a natural transformation $f^{[1]}:\mathcal U^{[1]}\Rightarrow\underline F$ \scth 
	\[
		f\circ\alpha^1-f\circ p_1=p_3\cdot f^{[1]}\ .
	\]
	We define the notion $\mathcal C^k_{MS}$ ($1\sle k\sle\infty$) inductively, just as above. A $\mathcal C^\infty_{MS}$ morphism is also called \emph{smooth}. By definition, any natural transformation $f:\mathcal U\Rightarrow\underline F$ in $\Top^\Lambda$ is $\mathcal C^0_{MS}$. 
	
	Let $k\in\nats\cup\{\infty\}$ and $1\sle\ell\sle k$, $\ell\in\nats$. For a $\mathcal C^k_{MS}$ morphism $f:\mathcal U\Rightarrow\underline F$, we define $d^\ell f:\mathcal U\times\underline E^\ell\to\underline F$ by $(d^\ell f)_\lambda=d^\ell f_\lambda$. Since $d^\ell f$ is given by a suitable restriction of $f^{[\ell]}:\mathcal U^{[\ell]}\to\underline F$, it is a $\mathcal C^{k-\ell}_{MS}$ morphism (where $\infty-\ell:=\infty$). The morphism $d^\ell f$ is called the $\ell$th derivative of $f$. 
\end{Defn}
	
\begin{Rem}
	The natural transformation $f^{[1]}$ is \emph{unique} whenever it exists; indeed, it is clear that if $f:\mathcal U\Rightarrow\underline F$ is $\mathcal C^1_{MS}$, then each $f_\lambda$ is $\mathcal C^1$ over $\lambda_0$, and one necessarily has that the $\lambda$-components of any $f^{[1]}$ are just the maps $f^{[1]}_\lambda$ occuring in the definition of `$\mathcal C^1$ over $\lambda_0$' for the maps $f_\lambda$ (where one manifestly has uniqueness).

	The above notions make sense if $f:\mathcal U\Rightarrow\underline F$ is considered as a morphism in the category $\Top^{\Lambda^\infty}$, if we assume that $R$, $E$ and $F$ are locally $k_\omega$. To this effect, we remark that finite products and countable inductive limits of locally $k_\omega$ spaces commute \cite[Proposition 4.7]{ggh-kac}. 
	
	We shall in the following not explicitly use the derivatives of $f$ introduced above, although some results could be stated with reference to them. We wish to emphasise, however, that to our opinion the relatively simple definition of the derivatives of a morphism is a particularily attractive feature of our version of the Molotkov--Sachse calculus. 
\end{Rem}

\begin{Prop}[c1natural]
	Let $E,F\in\TopSMod_R$, $\mathcal U\subset\underline E$ an open subfunctor, and $f:\mathcal U\Rightarrow\underline F$ a natural transformation. Then $f$ is $\mathcal C^k_{MS}$ if and only if \fa $\lambda\in\Lambda$, $f_\lambda$ is $\mathcal C^k$ over $\lambda_0$. If $R$, $E$ and $F$ are locally $k_\omega$, then equivalently, $f_\infty=f_{\lambda^\infty}$ is $\mathcal C^k$ over $\lambda_0^\infty$. 
\end{Prop}

\begin{proof}
	The only point one needs to check is that if $f_\lambda$ is $\mathcal C^1$ \fa $\lambda$, then $(f_\lambda^{[1]})$ is a natural transformation. This follows from \thmref{Lem}{c1natural} below.
\end{proof}

\begin{Par}
	Before we formulate the lemma, we need some terminology. 
	
	The object set of $\Lambda^\infty$ is ordered by inclusion. For $\lambda\in\Lambda^\infty$, let $\Lambda^\lambda$ denote the full subcategory of $\Lambda^\infty$ whose objects are those $\mu\in\Lambda^\infty$ for which $\mu\subset\lambda$. If $\lambda\in\Lambda$, then, of course, $\Lambda^\lambda\subset\Lambda$. 
	
	Let $E,F\in\TopSVec_R$, $\lambda\in\Lambda$ and $U\subset\underline E(\lambda)_{dW}$ be open. Let $\mathcal U$ be the associated functor in $\Top^\Lambda$, and let $\mathcal U_\lambda$, $\underline F_\lambda$ be the restrictions of $\mathcal U,\underline E$ to $\Lambda^\lambda$, considered as functors $\Lambda^\lambda\to\Sets$. We say that a map $f:U\to\underline F(\lambda)$ is \emph{natural} if it extends to a natural transformation $(f_\mu):\mathcal U_\lambda\to\underline F_\lambda$. 
	
		Let $\mu\subset\lambda$. Denote by $\eta_\mu:\mu\to\lambda$ the inclusion; since we have fixed the generators $\theta_i$, there exists a canonical epi $\eps_\mu:\lambda\to\mu$ whose kernel is the ideal of $\lambda$ generated by the $\theta_i$ which do not belong to $\mu$. \emph{E.g.}, $\eps_R=\eps$. 
\end{Par}

\begin{Lem}[c1natural]
	Let $U\subset\underline E(\lambda)$ be DeWitt-open, and let $f:U\to\underline F(\lambda)$ be $\mathcal C^1$ over $\lambda_0$ and natural. Then \fa $\mu\subset\lambda$, $f_\mu$ is $\mathcal C^1$ over $\mu_0$, and $f^{[1]}$ is natural. 
\end{Lem}

\begin{proof}
	We note that $\underline E(\eta_\mu)$ is simply the inclusion $\underline E(\mu)\subset\underline E(\lambda)$, so we will occasionally suppress it from the notation. It is clear by \thmref{Lem}{dewitt-top} that $U^{[1]}\subset\underline E(\lambda)^{[1]}$ (taken over $\lambda_0$) is DeWitt-open. Let $\mathcal U^{[1]}$ be the functor associated to $U^{[1]}$ by \thmref{Prop}{dewittgrothtop}. 
	
	Let $\mu\subset\lambda$. By naturality, 
	\[
		f_\mu(x+vt)-f_\mu(x)=f^{[1]}(x,v,t)\cdot t\mathfa (x,v,t)\in\mathcal U^{[1]}(\mu)\ .
	\]
	In particular, $f^{[1]}(U_\mu^{[1]})\subset\underline F(\mu)$, and we define $f_\mu^{[1]}=\underline F(\eps_\mu)\circ f^{[1]}\circ\underline E(\eta_\mu)$. Then $f_\mu$ and $f_\mu^{[1]}$ are continuous and satisfy (\ref{eq:diffdef}), so that $f_\mu$ is $\mathcal C^1$. 
	
	Let $\mu,\nu\subset\lambda$ and $\vphi:\mu\to\nu$ be a morphism. Write $\bar\vphi=\id\otimes\vphi$. We have 
	\begin{align*}
		f_\nu^{[1]}\Parens1{\bar\vphi(x),\bar\vphi(v),\vphi(t)}\cdot\vphi(t)&=f_\nu(\bar\vphi)(x+vt)-f_\nu(\bar\vphi(x))\\
		&=\bar\vphi\Parens1{f_\mu(x+vt)-f_\mu(x)}=\bar\vphi\Parens1{f^{[1]}_\mu(x,v,t)\cdot t}\\
		&=\bar\vphi\Parens1{f^{[1]}_\mu(x,v,t)}\cdot\vphi(t)
	\end{align*}
	Since $\mathcal U^{[1]}(\vphi)=\bar\vphi\times\bar\vphi\times\vphi$, this shows that 
	\[
		f_\nu^{[1]}\circ\mathcal U^{[1]}(\vphi)=\underline F(\vphi)\circ f_\mu^{[1]}\mathtxt{on}\mathcal U^{[1]}(\mu)\cap(\underline E(\mu)\times\underline E(\mu)\times\mu_0^\times)\ .
	\]
	Since $\mu_0^\times$ is dense, both sides of the equation are continuous, and the target spaces are Hausdorff, the equality holds everywhere, and $f^{[1]}$ is natural.
\end{proof}

	There is a strong version of the Taylor expansion for $\mathcal C^k_{MS}$ natural transformations. There is no remainder term if the order $k$ is high enough. Moreover, in the following proposition and corollary, we do \emph{not} have to assume that the integers $1,\dotsc,k$ are invertible in $R$. 

\begin{Prop}[taylor]
	Let $E,F\in\TopSMod_R$, $\mathcal U\subset\underline E$ an open subfunctor, and $f:\mathcal U\Rightarrow\underline F$. Let $f_\lambda$ be $\mathcal C^1$ over $R$ for every $\lambda\in\Lambda$. For any $\lambda\in\Lambda$, $x\in\mathcal U(\lambda)$, and $y\in\underline E(\lambda\theta_p)$ \fs $p$ \scth $\theta_p\in\lambda$, then 
	\[
		f_\lambda(x+y)=f_\lambda(x)+df_\lambda(x)y\ .
	\]
	If $R$, $E$ and $F$ are locally $k_\omega$, then it is sufficient to assume that $f_\infty=f_{\lambda^\infty}$ is $\mathcal C^1$ over $R$, and the conclusion holds for all $\lambda\in\Lambda^\infty$. 
\end{Prop}

\begin{proof}
	By naturality, $f_\lambda(x)\in\underline F(R[\theta_a|a\in A])$ whenever $A\subset\nats$ is a finite set \scth $x\in\mathcal U(R[\theta_a|a\in A])$ and $\lambda=\lambda^N$ for some $N\sge\max A$. 
	
	We first show that the $\theta_i$ occuring in $x$ and $y$, respectively, can be made independent of each other. In order to do so, we increase the number of variables and represent the $\theta_i$ occuring in $y$ as images of $\theta_i$ not occuring in $x$. Technically, the procedure goes as follows. 
	
	\Fs $N$, $\lambda=\lambda^N$, and $x\in\mathcal U(\lambda^N)$. We may write $y=\sum_{j=1}^My_j\theta_{I_j}\theta_p$ where $y_j\in E$ and $I_j=(i_{j1}<\dotsm<i_{jm_j})$. Let $m=\sum_{j=1}^M m_j$. We define the even unital algebra morphism $\alpha:\lambda^{N+m+1}\to\lambda^N$ which introduces our `new variables' as follows: Let $\alpha(\theta_i)=\theta_i$ \fa $i\sle N$, $\alpha(\theta_{N+m+1})=\theta_p$, and 
	\[
		\alpha\Parens1{\theta_{N+\sum_{j=1}^km_j+\ell}}=\theta_{i_{k+1,\ell}}\mathfa1\sle\ell\sle m_{k+1}\,,\,0\sle k<M\ .
	\]
	Set $z=y_1\theta_{N+1,\dotsc,N+m_1,N+m+1}+y_2\theta_{N+m_1+1,\dotsc,N+m_1+m_2,N+m+1}+\dotsm$. Then $\underline E(\alpha)(z)=y$ and $\mathcal U(\alpha)(x)=x$. Hence, we may first consider $z$ in place of $y$, and we will prove the main point of our assertion at this level.
	
	By naturality, applied to $\theta_{N+m+1}\mapsto0$, $\theta_\ell\mapsto\theta_\ell$ ($\ell\neq N+m+1$), 
	\[
	g(x,z)=f_{\lambda^{N+m+1}}(x+z)-f_{\lambda^{N+m+1}}(x)\in\underline F(\lambda^{N+m}\theta_{N+m+1})\ .
	\]
	Next, we scale $\theta_{N+m+1}$ with an element of the base ring: Let $c\in R$ and define $\beta:\lambda^{N+m+1}\to\lambda^{N+m+1}$ by $\beta(\theta_{N+m+1})=c\theta_{N+m+1}$ and $\beta(\theta_i)=\theta_i$ whenever $i\neq N+m+1$. Then, by naturality,
	\[
		cg(x,z)=\underline F(\beta)g(x,z)=f_{\lambda^{N+m+1}}(x+cz)-f_{\lambda^{N+m+1}}(x)=g(x,cz)\ ,
	\]
	and it follows that 
	\[
		cf^{[1]}_{\lambda^{N+m+1}}(x,z,c)=g(x,cz)=cg(x,z)=cf^{[1]}_{\lambda^{N+m+1}}(x,z,1)\ .
	\]
	Since $R^\times\subset R$ is dense, $f^{[1]}_{\lambda^{N+m+1}}(x,z,c)=f^{[1]}_{\lambda^{N+m+1}}(x,z,1)$. 
	
	In the last step, we apply $\alpha$ to return to our original setting with $y$ (instead of $z$). By naturality, 
	\[
		f^{[1]}_{\lambda^N}(x,y,c)=\underline F(\alpha)f^{[1]}_{\lambda^{N+m+1}}(x,z,c)=\underline F(\alpha)f^{[1]}_{\lambda^{N+m+1}}(x,z,1)=f^{[1]}_{\lambda^N}(x,y,1)\ .
	\]
	In particular, $f^{[1]}_{\lambda^N}(x,y,1)=f^{[1]}_{\lambda^N}(x,y,0)=df_{\lambda^N}(x)y$. This implies
	\[
		f_{\lambda^N}(x+y)-f_{\lambda^N}(x)=f^{[1]}_{\lambda^N}(x,y,1)=df_{\lambda^N}(x)y\ ,
	\]
	which is the assertion. 
	
	In the case that $R$, $E$ and $F$ are locally $k_\omega$, the final assertion follows by taking direct limits, in view of \cite[Corollary 5.7]{ggh-kac}. 
\end{proof}

\begin{Cor}[taylor]
	Let $f:\mathcal U\Rightarrow\underline F$ be given. If $f_\lambda$ is $\mathcal C^k$ over $R$ \fa $\lambda\in\Lambda$, where $k\in\nats$, then \fa $\lambda\in\Lambda$, $x\in\mathcal U(\lambda) $ and $y_j\in\underline E(\lambda\theta_{p_j})$, $\theta_{p_j}\in\lambda$, $j=1,\dotsc,k$, we have
	\[
		f_\lambda\Parens1{x+\textstyle\sum_{j=1}^ky_j}=f_\lambda(x)+\sum_{j=1}^k\sum_{\Abs0I=j}d^jf_\lambda(x)(y_{i_1},\dotsc,y_{i_j})
	\]
	where the second sum runs over all $I=(1\sle i_1<\dotsc<i_j\sle k)$. If $R$ and $E$ are locally $k_\omega$, then it is sufficient to assume that $f_\infty=f_{\lambda^\infty}$ is $\mathcal C^k$, and the conclusion holds for $\lambda\in\Lambda^\infty$. 
\end{Cor}

\begin{proof}
	This follows from \thmref{Prop}{taylor} by induction on $k$, since 
	\[
		f_\lambda\Parens1{x+\textstyle\sum_{j=1}^ky_j}=f_\lambda(x)+\sum_{j=1}^k\sum_{\Abs0I=j}df_\lambda^j(x)(y_{i_1},\dotsc,y_{i_j})
	\]
	implies, upon taking derivatives, 
	\[
		df_\lambda\Parens1{x+\textstyle\sum_{j=1}^ky_j}(y_{k+1})=df_\lambda(x)y_{k+1}+\sum_{j=1}^k\sum_{\Abs0I=j}df_\lambda^{j+1}(x)(y_{i_1},\dotsc,y_{i_j},y_{k+1})\ .
	\]
	Adding the equations gives the desired formula for $k+1$. 
\end{proof}

	A striking corollary of this exact Taylor expansion, in conjunction with the usual Taylor expansion, is that smoothness over $R$ implies smoothness over $\lambda_0$ ($\lambda\in\Lambda$) if the first derivative is already linear for the larger ring. For this, we need to assume that $R$ is a $\rats$-algebra. 

\begin{Prop}[smoothrlambda]
	Assume that $R$ is a $\rats$-algebra. Let $E,F\in\TopSMod_R$, $\mathcal U\subset\underline E$ be an open subfunctor, and $f:\mathcal U\Rightarrow\underline F$. If, \fa $\lambda\in\Lambda$, $f_\lambda$ is $\mathcal C^\infty$ over $R$ and the derivative $df_\lambda(x):\underline E(\lambda)\to\underline F(\lambda)$ is $\lambda_0$-linear \fa $x\in\mathcal U(\lambda)$, then \fa $\lambda\in\Lambda$, $f_\lambda$ is $\mathcal C^\infty$ over $\lambda_0$. If $R$, $E$ and $F$ are locally $k_\omega$, then it is sufficient to assume that $f_\infty=f_{\lambda^\infty}$ is $\mathcal C^\infty$ with $\lambda^\infty_0$-linear first derivatives at all points of $\mathcal U(\lambda^\infty)$, and the conclusion holds for all $\lambda\in\Lambda^\infty$. 
\end{Prop}

\begin{proof}
	First, we observe that all higher derivatives $d^kf_\lambda(x)(v_1,\dotsc,v_k)$ are $\lambda_0$-multilinear in $v_1,\dotsc,v_k$. Indeed, we certainly have $\lambda_0$-linearity in $v_1$ by assumption, and the higher derivatives are symmetric \cite[Lemma~4.8]{BGlN}. 
	
	We also record the following observation: If $t_1,\dotsc,t_{N+1}\in\lambda_0^N$ are nilpotent, then $t_1\dotsm t_{N+1}=0$ (since $\{1,\dotsc,N\}$ can be covered by at most $N$ disjoint non-void subsets). 
	
	Fix $\lambda\in\Lambda$. We write $U_\infty=\mathcal U(\lambda^\infty)$, $E_\infty=\underline E(\lambda^\infty)$, and $F_\infty=\underline F(\lambda^\infty)$. We claim now that there exist, for each positive integer $N$, continuous maps $g_N:U_N^{[1]}=\mathcal U^{[1]}(\lambda)\cap(U_\infty\times E_\infty\times\lambda^N_0)\to\underline F(\lambda)$ \scth
	\begin{equation}\label{eq:lambdanc1}
		f_\lambda(x+vt)-f_\lambda(x)=g_N(x,v,t)\cdot t\mathfa(x,v,t)\in U_N^{[1]}\ .
	\end{equation}
	
	Let $(x,v,t)\in U_N^{[1]}$ where $t=t_0+t^+$, $t_0\in R^\times$ and $t^+$ is nilpotent. By the usual Taylor expansion over $R$ \cite[Theorem~5.4]{BGlN}, there exists a map $R_N:U_0^{[1]}\to\underline F(\lambda)$ which is smooth over $R$, \scth
	\begin{align*}
		f_\lambda(x+vt)-f_\lambda(x)&=\begin{aligned}[t]\sum_{k=1}^{N+1}\frac{t_0^k}{k!}&d^kf_\lambda(x)(vt_0^{-1}t,\dotsc,vt_0^{-1}t)\\
		&+t_0^{N+1}R_N(x,vt_0^{-1}t,t_0)
		\end{aligned}\\
			&=\sum_{k=1}^{N+1}\frac{t^k}{k!}d^kf_\lambda(x)(v,\dotsc,v)+t_0^{N+1}R_N(x,vt_0^{-1}t,t_0)\ .
	\end{align*}
	(Note that $R_N$ is the $(N+1)$st Taylor remainder; and in the last equation, we have used the $\lambda_0$-multilinearity of the higher derivatives.)
	
	Because the remainder can be expressed via $f_\lambda$ and certain of its derivatives, the derivatives of $R_N(x,v,t)$ in the second argument are also $\lambda_0$-multilinear. 
	
	We may write $t^+=\sum_{j=1}^{N-1}t_j$ where $t_j\in R[\theta_k|k>j]\theta_j$; by \thmref{Cor}{taylor}, 
	\[
		t_0^{N+1}R_N(x,vt_0^{-1}t,t_0)=\begin{aligned}[t]&t_0^{N+1}R_N(x,v,t_0)\\
		&+\sum_{k=1}^{N-1}t_0^{N+1-k}\sum_{\Abs0I=k}t_Id^kR_N(x,\cdot,t_0)(v)(v,\dotsc,v)
		\end{aligned}
	\]
	where the latter sum extends over all multi-indices $I=(i_1<\dotsm<i_k)$. We have $t^{-1}=t_0^{-1}\sum_{j=0}^N(-1)^jt_0^{-j}(t^+)^j$, and for $I=(i_1<\dotsm<i_k)$, 
	\[
		t^{-1}t_0^{N+1-k}t_I=\sum_{j=0}^{N-k}(-1)^jt_0^{N-k-j}(t^+)^jt_I
	\]
	since $(t^+)^jt_I=0$ for $j+k>N$. 
	
	Therefore, we may define $g_{N+1}(x,v,t)=t^{-1}(f_\lambda(x+vt)-f_\lambda(x))$ for inver\-tible $t$, and this has a continuous extension to all of $U_N^{[1]}$. By \eqref{eq:lambdanc1}, $g_{N+1}$ extends $g_N$. By induction, $f_\lambda$ is $\mathcal C^1$ over $\lambda_0$.
	
	By its definition, it is clear that $g_N$ is $\mathcal C^\infty$ over $R$, and that the partial derivatives of $g_N(x,v,t)$ in $x$ and $v$ are $\lambda^N_0$-linear. Moreover, by the usual differentiation rules \cite[3.1, 3.3, 3.5]{BGlN}, we have 
	\[
		\tfrac d{dr}g_N(x,v,t+rs)|_{r=0}=t^{-1}df_\lambda(x+vt)(vs)-st^{-2}(f_\lambda(x+vt)-f_\lambda(x))
	\]
	for invertible $t$, and this is $\lambda^N_0$-linear in $s$ by the assumption on $f_\lambda$. It follows that $dg_N(x,v,t)$ is $\lambda^N_0$-linear, so that $f^{[1]}$ is $\mathcal C^\infty$ over $R$ with $\lambda_0$-linear derivative. The assertion now follows by a trivial induction. Finally, if $R$, $E$ and $F$ are locally $k_\omega$, then the above argument can be performed for $\lambda=\lambda^\infty$, and by taking direct limits. 
\end{proof}

\begin{Par*}
	We summarise our above considerations in the following theorem.
\end{Par*}

\begin{Th}[smoothequiv]
	Let $R$ be a $\rats$-algebra, $E,F\in\TopSMod_R$, and $\mathcal U\subset\underline E$ be an open subfunctor. Let maps $f_n:\mathcal U(\lambda^n)\to\underline F(\lambda^n)$ be given, \fa $n\in\nats$. The following are equivalent:
	\begin{enumerate}
		\item There exists a $\mathcal C^\infty_{MS}$ morphism $f:\mathcal U\Rightarrow\underline F$ in $\Top^\Lambda$, $f_{\lambda^n}=f_n$.
		\item The system $(f_n)$ defines a natural transformation from $\mathcal U$ to $\underline F$, and each of $f_n$ is $\mathcal C^\infty$ over $R$ with $\lambda^n_0$-linear derivatives at all points. 
	\end{enumerate}
	If $R$, $E$ and $F$ are locally $k_\omega$, then these statements are also equivalent to either of the following two:
	\begin{enumerate}
		\item[(i').] There exists a $\mathcal C^\infty_{MS}$ morphism $f:\mathcal U\Rightarrow\underline F$ in $\Top^{\Lambda^\infty}$, $f_{\lambda^n}=f_n$.
		\item[(ii').] There exists a natural map $f_\infty:\mathcal U(\lambda^\infty)\to\underline F(\lambda^\infty)$ extending the $f_n$ \scth $f_\infty$ is $\mathcal C^\infty$ over $R$ with $\lambda^\infty_0$-linear derivative at all points. 
	\end{enumerate}
\end{Th}

\begin{proof}
	Conditions (i) and (ii) are equivalent by \thmref{Prop}{smoothrlambda} and \thmref{Prop}{c1natural}. Assume that $R$, $E$ and $F$ are locally $k_\omega$. It is not hard to see, from the definitions and the fact that $\underline E(\lambda^\infty)$ is a topological $\lambda_0^\infty$-module, that one has (i) $\Leftrightarrow$ (i') and (ii) $\Leftrightarrow$ (ii'). The conclusion follows. 
\end{proof}

\begin{Rem}
	If $R=\reals$ or $R=\cplxs$, $E,F$ are $R$-Banach spaces, and $U\subset E$ is open, then a map $f:U\to F$ is $\mathcal C^k$ ($k\in\nats\cup\infty$) in the usual sense if and only if it is $\mathcal C^k$ in the sense of \thmref{Def}{c1def} \cite[Remark 7.3, Proposition 7.4]{BGlN}. Hence, \thmref{Th}{smoothequiv} shows that smoothness for natural transformations in the sense of \thmref{Def}{c1mordef} is equivalent to Sachse's definition of super-smoothness \cite[Definition 3.5.10]{sachse-diss}. 
	
	Let us point out the following slight defect of Sachse's definition: Given a morphism $f:\mathcal U\to\underline F$ which is super-smooth in the sense of Sachse, it is not clear from the definition that a super-smooth derivative $df:\mathcal U\times\underline E\to\underline F$ can be defined. However, this follows from \thmref{Th}{smoothequiv}: $df$ can be defined as a restriction of $f^{[1]}$, compare \thmref{Def}{c1mordef}. 
\end{Rem}

	We conclude the section with a more familiar formulation of the Taylor expansion, which will repeatedly be useful in the sequel.  

\begin{Prop}[taylor-ms]
	Let $R$ be a $\rats$-algebra, $E,F\in\TopSMod_R$, $\mathcal U\subset\underline E$ an open subfunctor and $f:\mathcal U\Rightarrow\underline F$ a smooth morphism. Fix elements $\lambda\in\Lambda$, $x\in\mathcal U(R)$, $n_0\in E_0\otimes\lambda_0^+$ and $n_1\in E_0\otimes\lambda_1$. Then
	\[
		f_\lambda(x+n_0+n_1)=\sum_{m,k=0}^\infty\frac1{m!\cdot k!}\cdot d^{m+k}f_\lambda(x)(n_0,\dotsc,n_0,n_1,\dotsc,n_1)\ ,
	\]
	where we evaluate $d^{m+k}f_\lambda(x)$ at $m$ copies of $n_0$ and at $k$ copies of $n_1$. 
\end{Prop}

\begin{proof}
	The point to note is that if $b$ is a symmetric $\rats$-$k$-linear map, then
	\[
		\sum_{1\sle i_1\sle\dotsm\sle i_k\sle n}b(x_{i_1},\dotsc,x_{i_k})=\frac1{k!}\cdot b\Parens1{\textstyle\sum_{j=1}^nx_j,\dotsc,\sum_{j=1}^nx_j}\ .
	\]
	Moreover, in our application, the ascending multi-indices which are not strictly ascending do not contribute. Then the formula follows immediately from \thmref{Cor}{taylor}.
\end{proof}

\begin{Rem}
	The formula in \thmref{Prop}{taylor-ms} is used without proof in \cite[Theorem~4.11]{sachse-preprint} and \cite[(3.3.1)]{molotkov}.
\end{Rem}

\section{Equivalence of categories of supermanifolds}

From hereon, we assume that $R$ is a unital commutative Hausdorff topo\-lo\-gi\-cal $\rats$-algebra with a dense group of units. 

\subsection{Equivalence of categories of superdomains}\label{section31}

\begin{Defn}[Categories]
We define the category $\SDom_{MS}=\SDom_{MS}(R)$ of superdomains (over $R$) in the sense of Molotkov--Sachse, as follows: The objects are pairs $(\mathcal U,E)$ where $E\in\TopSMod_R$ and $\mathcal U\subset\underline E\in\Top^\Lambda$ is an open subfunctor; a morphism $f\colon(\mathcal U,E)\to(\mathcal V,F)$ is a natural transformation $f\colon\mathcal U\Rightarrow\mathcal V$ which is smooth (\emph{i.e.}~$\mathcal C^\infty_{MS}$) when considered as a natural transformation $f\colon\mathcal U\Rightarrow\underline F$. The full subcategory $\SDom_{MS}^{fd}$ of \emph{finite-dimensional} superdomains consists of those pairs $(\mathcal U,E)$ where $E$ is finite-dimensional. 

We shall usually suppress the mention of $E$ in our notation. Indeed, since $\mathcal U(\lambda^2)=\mathcal U(R)+E_0\theta_1\theta_2\oplus E_1\theta_1\oplus E_1\theta_2$, $E$ is uniquely determined by $\mathcal U$ as a topological space. In particular, we write $\Ct[^\infty_{MS}]0{\mathcal U,\mathcal V}$ for the set of morphisms $(\mathcal U,E)\to(\mathcal V,F)$. 

If $R$ is a field, then the category $\SDom_{BKL}=\SDom_{BKL}(R)$ of superdomains (over $R$) in the sense of Berezin--Kostant--Leites is given as follows: Objects are pairs $(U,\mathcal F)$ where $U$ is an open subset of some \emph{finite-dimensional} $R$-vector space $E_0$, and $\mathcal F$ is a sheaf of supercommutative superalgebras isomorphic to $\mathcal C^\infty_U\otimes\bigwedge E_1^*$ where $E_1$ is another finite-dimensional $R$-vector space (the sheaf property follows from  \cite[Lemma 4.9]{BGlN}); morphisms are pairs $(f,f^*):(U,\mathcal F)\to(V,\mathcal G)$ where $f:U\to V$ is continuous and $f^*:\mathcal G\to f_*\mathcal F$ is an even morphism of sheaves of unital superalgebras. 
\end{Defn}

\begin{Rem}
If $h\in\mathcal F_x$ where $x\in U$ and $\mathcal F_x$ is the stalk at $x$, then $h$ is invertible if and only if $h\not\in\ger m_{\mathcal F,x}$ where the latter denotes the maximal ideal. It follows that any morphism in $\SDom_{BKL}$ is \emph{local}, \emph{i.e.}~$f^*(\ger m_{\mathcal F,x})\subset\ger m_{\mathcal G,f(x)}$. 
\end{Rem}

In the course of the present section we shall show that the categories $\SDom_{MS}^{fd}$ and $\SDom_{BKL}$ are equivalent, by explicitly defining (in \thmref{Def}{functorphidef}) a functor $\Phi:\SDom_{MS}^{fd}\to\SDom_{BKL}$, and proving that it is an equivalence. To this end, a crucial step will be an alternative description of the morphisms in the Molotkov--Sachse category, which we presently derive. 

\begin{Par}
	Let $\mathcal U,\mathcal V\in\SDom_{MS}$ and $W\subset U=\mathcal U(R)$ be open. We set 
	\[
		\mathcal O_{\mathcal U,\mathcal V}(W)=\Set1{(\vphi_k)_{k\sge0}}{\vphi_0:W\to V\,,\,\vphi_k:W\to\Alt[^k]0{E_1,F_k}\ \text{smooth}}\ .
	\]
	Here, $\mathcal U\subset\underline E$, $\mathcal V\subset\underline F$ are open subfunctors, $V=\mathcal V(R)$, $\Alt[^k]0{E_1,F_k}$ denotes the set of alternating multilinear maps $E_1^k\to F_k$ ($F_k=F_0$ or $F_k=F_1$ according to the parity of $k$), and $\vphi_k:W\to\Alt[^k]0{E_1,F_k}$ is called smooth if so is the map $W\times E_1^k\to F_k:(x,v_1,\dotsc,v_k)\mapsto\vphi_k(x)(v_1,\dotsc,v_k)$. 

	Let $\mathcal R\in\Top^{\Lambda^\infty}$ be defined by $\mathcal R(\lambda)=\lambda$ and $\mathcal R(\alpha)=\alpha$ for any arrow $\alpha$ in $\Lambda^\infty$. Then $\mathcal R$ is naturally equivalent to $\underline{R^{1|1}}$, and we may define as a particular case $\mathcal O_{\mathcal U}(W)=\mathcal O_{\mathcal U,\mathcal R}(W)$. Restriction maps on $\mathcal O_{\mathcal U,\mathcal V}$ are defined in the obvious way, and this gives rise to a presheaf on $U$. If $R$ is a \emph{field} (of characteristic zero), then by \cite[Lemma 4.9]{BGlN}, $\mathcal O_{\mathcal U,\mathcal V}$ is a sheaf.

	On $\mathcal O_{\mathcal U}(W)$, we may define an algebra structure pointwise. \emph{I.e.}, $\ger S_{k,\ell}$ denoting $(k,\ell)$ shuffles, 
	\[
		(\vphi\cdot\psi)_m(x,v_1,\dotsc,v_m)=\sum_{k=0}^m\sum_{\sigma\in\ger S_{k,m-k}}
		\begin{aligned}[t]
			\eps(\sigma)&\cdot\vphi_k(x)(v_{\sigma(1)},\dotsc,v_{\sigma(k)})\\
			&\cdot \psi_{m-k}(x)(v_{\sigma(k+1)},\dotsc,v_{\sigma(m)})\ .
		\end{aligned}
	\]
	
	If we define a $\ints_2$-grading \emph{via}
	\[
	\vphi\in\mathcal O_{\mathcal U}(W)_i\quad\Leftrightarrow\quad\vphi_k=0\quad(\forall k\ \text\scth k+i\equiv1\ (2))\ ,
	\] 
	then $\mathcal O_{\mathcal U}(W)$ is a supercommutative superalgebra. Let $\mathcal C^\infty_U$ be the subpresheaf of $\mathcal O_{\mathcal U}$ which consists of all $\vphi=(\vphi_k)$ \scth $\vphi_k=0$ for $k>0$. It is a commutative subalgebra. 
	
	Furthermore, for any open $W\subset U$, we let 
	\[
		\phi=\phi_{\mathcal U}(W):\Ct[^\infty_{MS}]0{\mathcal U_W,\mathcal V}\to\mathcal O_{\mathcal U,\mathcal V}(W)
	\]
	be defined by the equation $\phi(f)=\vphi$, 
	\[
		f_{\lambda^\infty}\Parens1{x+\textstyle\sum_{j=1}^Ny_j\theta_j}=\sum_{k=0}^N\sum_{\Abs0I=k}\vphi_k(x)(y_{i_1},\dotsc,y_{i_k})\theta_I
	\]
	\fa $N\in\nats$, $x\in U$, $y_j\in E_1$, where as usual, the inner sum extends over all $I=(1\sle i_1<\dotsm<i_k\sle N)$. 
\end{Par}

\begin{Prop}[structuresheafiso]
	Let $\mathcal U,\mathcal V\in\SDom_{MS}$. For any open $W\subset U$, the map 
	\[
		\phi=\phi_{\mathcal U}(W):\Ct[^\infty_{MS}]0{\mathcal U_W,\mathcal V}\to\mathcal O_{\mathcal U,\mathcal V}(W)
	\]
	is a bijection. 
\end{Prop}

\begin{proof}
	For $\vphi\in\mathcal O_{\mathcal U}(W)$, define $f_\vphi\colon\mathcal U_W(\lambda^\infty)\to\mathcal V(\lambda^\infty)$ by
	\begin{equation}\label{eq:grassmannanalytic}
		f_{\vphi}\Parens1{x+n_0+n_1}=\sum_{m,k=0}^\infty\frac1{m!\cdot k!}\cdot d^m\vphi_k(x)(n_0,\dotsc,n_0,n_1,\dotsc,n_1)
	\end{equation}
	\fa $x\in W$, $n_0\in E_0\times\lambda_0^{\infty+}$, $n_1\in E_1\otimes\lambda_1^\infty$. Here, it is understood that 
	\[
	d^m\vphi_k(x)(v_1\theta_{I_1},\dotsc,v_{k+m}\theta_{I_{k+m}})=d^m\vphi_k(x)(v_1,\dotsc,v_{k+m})\theta_{I_1}\dotsm\theta_{I_{k+m}}
	\]
	\fa $v_1,\dotsc,v_k\in E_0$, $v_{k+1},\dotsc,v_{k+m}\in E_1$, where $\Abs0{I_j}\equiv0\ (2)$ for $j\sle k$ and $\Abs0{I_j}\equiv1\ (2)$ for $j>k$. (This makes $d^m\vphi_k(x)(n_0,\dotsc,n_0,n_1,\dotsc,n_1)$ \emph{symmetric} both in the $n_0$ and in the $n_1$ variables.)
	
	We claim that $f_\vphi$ defines a smooth morphism $\mathcal U_W\Rightarrow\mathcal V$. By its definition, it is straightforward to see that $f_\vphi$ defines, by restriction onto $\mathcal U_W(\lambda)$, where $\lambda\in\Lambda$, a natural transformation $\mathcal U_W\Rightarrow\mathcal V$, if these are considered as set-valued functors $\Lambda\to\mathrm{Sets}$. By \thmref{Th}{smoothequiv}, the remaining issue is to check whether $f_\vphi|_{\mathcal U_W(\lambda)}$ is smooth over $R$ with a $\lambda_0$-linear derivative at any point. It is clear that $f_\vphi|_{\mathcal U_W(\lambda)}$ is smooth, since so are all the $\vphi_k$.  
	
	We compute the derivative. Let $x\in W$, $v\in E_0$, $x_i,v_i\in E_i\otimes\lambda_i^+$. Set $y=x+x_0+x_1$, $u=v+v_0+v_1$. Then 
	\[
		df_\vphi(y)(u)=\sum_{m,k=0}^\infty\frac1{m!\cdot k!}\cdot
		\begin{aligned}[t]
			\bigl[&d^{m+1}\vphi_k(x)(v,x_0,\dotsc,x_0,x_1,\dotsc,x_1)\\
			&+m\cdot d^m\vphi_k(x)(v_0,x_0,\dotsc,x_0,x_1,\dotsc,x_1)\\
			&+k\cdot d^m\vphi_k(x)(x_0,\dotsc,x_0,v_1,x_1,\dotsc,x_1)\bigr]\ .
		\end{aligned}
	\]
	Here, it is understood that the last two summands are zero for $m=0$ and $k=0$, respectively. 
	
	This expression is certainly $\lambda_0$-linear. \emph{E.g.}, for any $a\in\lambda_0^+$, we have $av\in E_0\otimes\lambda^+_0$, and 
	\begin{multline*}
		a\cdot\frac1{m!\cdot k!}\cdot d^{m+1}\vphi_k(x)(v,x_0,\dotsc,x_0,x_1,\dotsc,x_1)\\
		=\frac{m+1}{(m+1)!\cdot k!}\cdot d^{m+1}\vphi_k(x)(av,x_0,\dotsc,x_0,x_1,\dotsc,x_1)\ ,
	\end{multline*}
	where the right hand side occurs in the expression for $df_\vphi(y)(au)$. Thus, $f_\vphi$ defines a smooth morphism $\mathcal U_W\Rightarrow\mathcal R$. 
	
	By construction, it is immediate that $\vphi=\phi(f_\vphi)$. Conversely, using \eqref{eq:grassmannanalytic}, \thmref{Th}{smoothequiv} and  \thmref{Prop}{taylor-ms}, it follows that if $\vphi=\phi(f)$ \fs morphism $f\in\Ct[^\infty_{MS}]0{\mathcal U_W,\mathcal V}$, then $f=f_\vphi$. Hence, $\phi$ is an isomorphism.
\end{proof}

\begin{Rem}
	Equation \eqref{eq:grassmannanalytic} defines exactly the well-known `Grassmann-analytic continuation' due to Berezin \cite{berezin}. 
\end{Rem}

\begin{Cor}
	The map $\phi:\Ct[^\infty_{MS}]0{\mathcal U_W,\mathcal R}\to\mathcal O_{\mathcal U}(W)$ is an isomorphism of unital $R$-algebras. Here, the product on $\Ct[^\infty_{MS}]0{\mathcal U_W,\mathcal R}$ is defined by
	\[
	(f\cdot g)_\lambda(x)=f_\lambda(x)g_\lambda(x)\mathfa f,g\in\Ct[^\infty_{MS}]0{\mathcal U_W,\mathcal R}\,,\,\lambda\in\Lambda\,,\,x\in\mathcal U(\lambda)\ .
	\]
\end{Cor}

\begin{proof}
	We need only check that it is an algebra morphism. Let $\vphi=\phi(f)$, $\psi=\phi(g)$, and $z=x+\sum_{j=1}^Ny_j\theta_j$. We abbreviate 
	\[
	f_k(z)=\sum_{\Abs0I=k}\vphi_k(x)(y_{i_1},\dotsc,y_{i_k})\theta_I
	\] 
	where the sum ranges over all $I=(1\sle i_1<\dotsm<i_k\sle N)$. Then 
	\begin{align*}
		(f\cdot g)_{\lambda^N}(z)&=\sum_{m=0}^{2N}\sum_{k+\ell=m}f_k(z)g_\ell(z)\\
		&=\sum_{m=0}^{2N}\sum_{\Abs0I=m}\sum_{k+\ell=m}\sum_{\sigma\in\ger S_{k,\ell}}
			\begin{aligned}[t]
				\eps(\sigma)&\cdot\vphi_k(x)(y_{i_{\sigma(1)}},\dots,y_{i_{\sigma(k)}})\\
				&\cdot\psi_\ell(x)(y_{i_{\sigma(k+1)}},\dotsc,y_{i_{\sigma(m)}})\theta_I
			\end{aligned}\\
		&=\sum_{m=0}^N\sum_{\Abs0I=m}(\vphi\cdot\psi)_m(x)(y_{i_1},\dotsc,y_{i_m})\theta_I
	\end{align*}
	since $\theta_I=0$ for $m>N$. This proves the claim. 
\end{proof}

Using the above description of the morphisms in the Molotkov--Sachse category, one can derive a formula for the composition of two morphisms. 

\begin{Prop}[composfmla]
	Let $\mathcal U,\mathcal V,\mathcal W\in\SDom_{MS}$, $f\in\Ct[^\infty_{MS}]0{\mathcal U,\mathcal V}$, $g\in\Ct[^\infty_{MS}]0{\mathcal V,\mathcal W}$. Define $\vphi=\phi(f)$, $\psi=\phi(g)$ and $\vrho=\phi(g\circ f)$. 
	Then \fa $n\in\nats$, and all $v=(v_1,\dotsc,v_n)\in E_1^n$, 
	\begin{equation}\label{eq:skelcompos}
		\vrho_n(x)(v)=\sum_{\substack{m,k,(\alpha,\beta)\in I_{m,k}^n,\\\sigma\in\ger S_{\Abs0\alpha},\tau\in\ger S_{\Abs0\beta}}}\frac{\eps(\tau)}{m!k!\alpha!\beta!}d^m\psi_k(\vphi_0(x))\Parens1{(\vphi_\alpha\times\vphi_\beta)(x)(v^{\sigma,\tau})}
	\end{equation}
	where
	\begin{gather*}
		I_{m,k}^n=\Set1{(\alpha,\beta)\in(2\nats)^m\times(2\nats+1)^k}{\forall j\,:\,\alpha_j>0\ ,\ \Abs0\alpha+\Abs0\beta=n}\ ,\\
		\vphi_\alpha=\vphi_{\alpha_1}\times\dotsm\times\vphi_{\alpha_n}\ ,\\
		v^{\sigma,\tau}=(v_{\sigma(1)},\dotsc,v_{\sigma(\Abs0\alpha)},v_{\Abs0\alpha+\tau(1)},\dotsc,v_{\Abs0\alpha+\tau(\Abs0\beta)})\ .
	\end{gather*}
\end{Prop}

\begin{proof}
	As we know, $\vrho$ is determined by the equation
	\[
		g_\lambda(f_\lambda(x+y))=\sum_{k=0}^\infty\frac1{k!}\vrho_k(x)(y,\dotsc,y)
	\]
	where $\lambda=\lambda^n$, $x\in U=\mathcal U(R)$, and $y=\sum_{j=1}^ny_j\theta_j\in E_1\otimes(\lambda^n)_1$. Let $n_i=\sum_{0\neq k\in2\nats+i}\frac1{k!}\vphi_k(x)(y,\dotsc,y)$. Then 
	\[
		g_\lambda(f_\lambda(x+y))=\sum_{m,k=0}^\infty\frac1{m!k!}d^m\psi_k(\vphi_0(x))(n_0,\dotsc,n_0,n_1,\dotsc,n_1)\ .
	\]
	From this, it is not hard to deduce the formula. 
\end{proof}

\begin{Rem}
	The formula \eqref{eq:skelcompos} in \thmref{Prop}{composfmla} is stated without proof in \cite[Proposition 3.3.3]{molotkov}. 
\end{Rem}

So far, we have derived an alternative description of the morphisms in the Molotkov--Sachse category. This will enable us to identify morphisms of \emph{finite-dimensional} Molotkov--Sachse superdomains and of Berezin--Leites superdomains. The following proposition defines the sought-for equivalence $\Phi:\SDom_{MS}^{fd}\to\SDom_{BKL}$ on the level of objects; moreover, it indicates that $\Phi$ is an adjoint of the point functor explained in the introduction.  

\begin{Prop}[lambdapoints]
	Assume that $R$ is a field. Let $\lambda\in\Lambda$ and define the superdomain $*_\lambda=(*,\lambda)\in\SDom_{BKL}$ whose underlying domain is a point and whose structure sheaf is the constant sheaf $\lambda$.  
	
	For any $\mathcal U\in\SDom_{MS}^{fd}$, let $\Phi(\mathcal U)=(U,\mathcal O_\mathcal U)$ where $U=\mathcal U(R)$. The map
	\[
		\eps:\mathcal U(\lambda)\to\Hom[_{BKL}]0{*_\lambda,\Phi(\mathcal U)}
	\]
	where for $x=x_R+x^+$, we let
	\[
	\eps_x=(\eps_{x_R},\eps_x^*)\,,\,\eps_{x_R}(*)=x_R\,,\,\eps_x^*\phi(f)=f_\lambda(x)\ ,
	\]
	is a bijection. 
\end{Prop}

	In the \emph{proof}, we need two lemmata. Both do not generalise to infinite dimensions. This is the main reason why the ringed space approach to (super-) manifolds is not well-behaved in this case. 

\begin{Lem}[hadamard]
	Let $E_0$ be a finite-dimensional $R$-vector space, $U\subset E_0$ be open and $x\in U$. Let $\mathcal C^\infty_{U,x}$ be the stalk of $\mathcal C^\infty_U$ at $x$. Then the  maximal ideal $\ger m_x=\Set0{f\in\mathcal C^\infty_{\mathcal U,x}}{f(x)=0}$ is generated by $\mu-\mu(x)$ where $\mu\in E_0^*$.
\end{Lem}

\begin{proof}
	Certainly, $\mu-\mu(x)\in\ger m_x$ for $\mu\in E_0^*$. Conversely, let $f\in\ger m_x$. This germ is represented by a smooth function defined on an open neighbourhood $W\subset U$ of $x$. Let $e_1,\dotsc,e_n$ be a basis of $E_0$. For $y\in W$, we have 
	\begin{align*}
		f(y)&=f\Parens1{x+\textstyle\sum_{j=1}^n(y_j-x_j)e_j}\\
		&=
			\begin{aligned}[t]
				(y_1-x_1)&\cdot f^{[1]}\Parens1{x+\textstyle\sum_{j=2}^n(y_j-x_j)e_j,e_1,y_1-x_1}\\
				&+f\Parens1{x+\textstyle\sum_{j=2}^n(y_j-x_j)e_j}\ .
			\end{aligned}
	\end{align*}
	Let $h_k(y)=f^{[1]}\Parens1{x+\textstyle\sum_{j=k+1}^n(y_j-x_j)e_j,e_k,y_k-x_k}$. By induction,
	\[
		f(y)=\sum_{j=1}^n(y_j-x_j)h_j(y)+f(x)=\sum_{j=1}^n(y_j-x_j)h_j(x)\ .
	\]
	This proves the assertion.
\end{proof}

\begin{Lem}[smoothfnlambdamor]
	Let $E_0$ be a finite-dimensional $R$-vector space, $U\subset E_0$ be open, $x_0\in U$, and $\lambda\in\Lambda$. With any algebra morphism $a\colon\mathcal C^\infty_{U,x_0}\to\lambda$, there is associated a unique $x\in x_0+E_0\otimes\lambda^+$; if $a$ is even, then $x\in x_0+E_0\otimes\lambda_0^+$. The correspondence is given by 
	\[
		x=\textstyle\sum_Ix_I\theta_I\ ,\ a(\mu)=\sum_I\mu(x_I)\theta_I\mathfa\mu\in E_0^*\ .
	\]
	Moreover, $a$ is uniquely determined by $x$. 
\end{Lem}

\begin{proof}
	Let $a:\mathcal C^\infty_{U,x_0}\to\lambda$ be an algebra morphism. Observe that we have $\lambda=R[\theta_1,\dotsc,\theta_N]$ \fs $N$. We decompose $a(f)=\sum_Ia_I(f)\theta_I$ where $a_I(f)\in R$. Then for $m=\Abs0I$
	\begin{equation}\label{eq:algmorrecursion}
		a_I(fg)=\sum_{k=0}^m\sum_{\sigma\in\ger S_{k,m-k}}\eps(\sigma)a_{i_{\sigma(1)}\dotsm i_{\sigma(k)}}(f)a_{i_{\sigma(k+1)}\dotsm i_{\sigma(m)}}(g)\ .
	\end{equation}
	In particular, $a_0$ is an algebra morphism, and $a_1(1)=2a_1(1)=0$. Inductively, let $\Abs0I>0$ and assume that $a_J(1)=0$ for all $\Abs0I>\Abs0J>0$. Then 
	\[
		a_I(1)=a_0(1)a_I(1)+a_I(1)a_0(1)=2a_I(1)=0\ .
	\]
	Hence, $a_I$ ($\Abs0I>0$) is determined by its values on $\ger m_{x_0}$. The germ $f\in\mathcal C^\infty_{U,x_0}$ is invertible if and only if $f(x_0)\neq0$. Thus, $a_0(f)-f(x_0)=a_0(f-f(x_0))=0$. 
	
	Since $\ger m_{x_0}$ is generated as an ideal by $\mu-\mu(x_0)$, $\mu\in E_0^*$, and we have 
	\[
		a_1((\mu-\mu(x_0))h)=a_0(\mu-\mu(x_0))a_1(h)+a_1(\mu-\mu(x))a_0(h)=a_1(\mu)h(x_0)\ ,
	\]
	it follows that $a_1$ is determined by its values on $E_0^*$. An easy induction using \eqref{eq:algmorrecursion} shows that $a_I$, $\Abs0I>0$, is determined by its values on $E_0^*$. But since $E_0^{**}=E_0$, there are unique $x_I\in E_0$ \scth $a_I(\mu)=\mu(x_I)$ \fa $\mu\in E_0^*$. Thus, $a$ is uniquely determined by $x=\sum_Ix_I\theta_I\in x_0+E_0\otimes\lambda^+$.
\end{proof}

\begin{proof}[Proof of \thmref{Prop}{lambdapoints}]
	The map is certainly well-defined. Let $(\alpha,\alpha^*)$ be a morphism $*_\lambda\to\Phi(\mathcal U)$. Then we have a point $x_0=\alpha(*)\in U$. 
	
	The value $\alpha^*(\vphi)$ (for $\vphi\in\mathcal O_{\mathcal U}(W)$, $x_0\in W\subset U$) depends only on the germ of $\vphi$ at $x_0$ since for any open $W'\subset W$, $x_0\in W'$, the following diagram commutes,
	\[
		\xymatrix@C+5ex{\mathcal O_{\mathcal U}(W)\ar[r]^-{\alpha^*_W}\ar[d]_-{\vrho_{W'}^W}&\lambda\\
		\mathcal O_{\mathcal U}(W')\ar[ru]_-{\alpha^*_{W'}}}
	\]
	(because $\alpha^*:\mathcal O_{\mathcal U}\to\alpha_*\lambda$ is a sheaf morphism). 
	
	Now, $\mathcal O_{\mathcal U,x_0}=\mathcal C^\infty_{U,x_0}\otimes\bigwedge E_1^*$. By \thmref{Lem}{smoothfnlambdamor}, there is a unique element $x_0^+\in E_0\otimes\lambda_0^+$ \scth $\alpha^*\vphi=f_\lambda(x_0+x_0^+)$ \fa $\vphi=\phi(f)\in\mathcal C^\infty_{\mathcal U,x_0}$. On the other hand, the set of even algebra morphisms $\bigwedge E_1^*\to\lambda$ equals 
	\[
		\Hom0{E_1^*,\lambda}_0=(E_1\otimes\lambda)_0=E_1\otimes\lambda_1\ .
	\]
	Thus, there is a unique element $x_1^+=\sum_{\Abs0I\equiv1\ (2)}x_I\theta_I\in E_1\otimes\lambda_1$ \scth \fa $\vphi=\phi(f)\in\bigwedge E_1^*$, $\alpha^*\vphi=f_\lambda(x_1^+)$. 

	If $\vphi=\phi(f)\in\mathcal C^\infty_{U,x}$ and $\psi=\phi(g)\in\bigwedge E_1^*$, then $f_\lambda(x_0+x_0^++x_1^+)=f_\lambda(x_0+x_0^+)$ and $g_\lambda(x_0+x_0^++x_1^+)=g_\lambda(x_1^+)$, as follows from \eqref{eq:grassmannanalytic}. Thus, 
	\[
	\alpha^*(\vphi\psi)=\alpha^*\vphi\cdot\alpha^*\psi=f_\lambda(x_0+x_0^+)g_\lambda(x_1^+)=(f\cdot g)_\lambda(x_0+x_0^++x_1^+)\ .
	\]
	Since $\vphi\psi=\phi(fg)$, it follows that $\alpha^*=\eps_x^*$ where $x=x_0+x_0^++x_1^+$. Since $x$ was by construction unique, it follows that $\eps$ is a bijection. 
\end{proof}

Finally, we are in a position to define the functor $\Phi$, and to prove that it is an equivalence of categories. 
	
\begin{Defn}[functorphidef]
	Assume that $R$ is a field. Define a functor $\Phi$ from $\SDom_{MS}^{fd}$ to $\SDom_{BKL}$ as follows: On objects, we let $\Phi(\mathcal U)=(U,\mathcal O_{\mathcal U})$. Any $f\in\Ct[^\infty_{MS}]0{\mathcal U,\mathcal V}$ is mapped to the morphism $\Phi(f)=(f_R,f^*):\Phi(\mathcal U)\to\Phi(\mathcal V)$ given by  
	\[
		f^*\phi(h)=\phi\Parens1{h\circ f|_{\mathcal U_{f^{-1}(W)}}}\mathfa h\in\Ct[^\infty_{MS}]0{\mathcal V_W,\mathcal R}\ ,
	\]
	where $W$ runs through the open subsets of $V=\mathcal V(R)$. 
\end{Defn}

\begin{Th}[sdomequiv]
	Assume that $R$ is a non-discrete Hausdorff topological field of characteristic zero. Then $\Phi:\SDom^{fd}_{MS}\to\SDom_{BKL}$ is an equivalence of categories. 
\end{Th}

In the \emph{proof}, we need to extend Hadamard's lemma (\emph{cf.}~\cite{leites,schmitt-supergeom}) in the usual way.

\begin{Lem}[genhadamard]
	Let $\underline E\supset\mathcal U\in\SDom_{MS}^{fd}$, $U=\mathcal U(R)$, and $x\in U$. Let $\ger m_{\mathcal O,x}$ denote the maximal ideal of $\mathcal O_{\mathcal U,x}$. 
	\begin{enumerate}
		\item As an ideal, $\ger m_{\mathcal O,x}$ is generated by $\mu-\mu(x)$, $\nu$ for $\mu\in E_0^*$, $\nu\in E_1^*$.
		\item Let $W\subset U$ be an open neighbourhood of $x$. For any $f\in\mathcal O_{\mathcal U}(W)$ and any $n$, there exists a polynomial $p$ of degree $\sle n$ in $\mu-\mu(x)$, $\nu$ (with coefficients in $R$), where $\mu\in E_0^*$, $\nu\in E_1^*$, \scth one has $(f-p)_x\in\ger m_{\mathcal O,x}^{n+1}$ where $f_x$ denotes the germ at $x$. 
		\item Let $W\subset U$ be open and $f\in\mathcal O_{\mathcal U}(W)$. If $q=\dim E_1$ and $f_y\in\ger m_{\mathcal O,x}^{q+1}$ \fa $y\in W$, then $f=0$. 
	\end{enumerate}
\end{Lem}

\begin{proof}
	Statement (i) follows easily from \thmref{Lem}{hadamard}, and (ii), (iii) can then be deduced by standard procedures \cite[Lemma 2.14 and proof]{schmitt-supergeom}.
\end{proof}

\begin{proof}[Proof of \thmref{Th}{sdomequiv}]
	It is clear that $\Phi$ is an essentially surjective functor. If $f\colon\mathcal U\Rightarrow\mathcal V$ is a morphism in $\SDom^{fd}_{MS}$, then $f$ is determined uniquely by the maps 
	\[
	\Hom0{*_\lambda,\Phi(\mathcal U)}\to\Hom0{*_\lambda,\Phi(\mathcal V)}\colon\psi\mapsto\Phi(f)\circ\psi\ , 	\]
	for $\lambda\in\Lambda$, by \thmref{Prop}{lambdapoints}. This proves that $\Phi$ is faithful. 
	
	It remains to show that $\Phi$ is full. To that end, let $(\psi,\psi^*):\Phi(\mathcal U)\to\Phi(\mathcal V)$ be a morphism in $\SDom_{BKL}$. We let $U=\mathcal U(R)$ and $V=\mathcal V(R)$, and define $\vphi_0=\psi:U\to V$. As in the proof of \thmref{Lem}{smoothfnlambdamor}, it follows that $h\circ\vphi_0=(\psi^*h)_0$ \fa $h\in\mathcal C_V^\infty\subset\mathcal O_{\mathcal V}$. In particular, $\mu\circ\vphi_0:U\to R$ is a smooth function \fa $\mu\in F_0^*$. Since any map $g:U\to F_0$ is continuous if and only if $\mu\circ g$ is continuous \fa $\mu\in F_0^*$ ($F_0$ carries the product topology with respect to any chosen linear isomorphim with $R^n$, $n=\dim F_0$), it follows easily that $\vphi_0$ is smooth. 
	
	Next, a simple induction shows that for all \emph{even} $n>0$ there are unique smooth maps $\vphi_n:U\to\Alt[^n]0{E_1,F_0}$ \scth \fa $\mu\in F_0^*$, $x\in U$, $v=(v_1,\dotsc,v_n)\in E_1^n$, 
	\[
		\mu(\vphi_n(x)(v))=[\psi^*\mu(x)]_n(v)-\mu(\vphi_0(x))-\sum_{\substack{2\sle m\sle n/2,\\\alpha\in I^n_{m,0},\sigma\in\ger S_n}}\frac1{m!\alpha!}\mu(\vphi_\alpha(x)(v^\sigma))\ .
	\]
	Here, we identify $\Ct[^\infty]0{U,\bigwedge E_1^*}\cong\Ct[^\infty]0U\otimes\bigwedge E_1^*=\mathcal O_{\mathcal U}(U)$, and $[\cdot]_n$ denotes the homogeneous component of degree $n$. Moreover, we use the notation from \thmref{Prop}{composfmla}. We observe that on the right hand side, $\vphi_\alpha$ has only components $\vphi_{\alpha_j}$ where $\alpha_j<n$. 
	
	For all \emph{odd} $n$, there are unique smooth maps $\vphi_n:U\to\Alt[^n]0{E_1,F_1}$ \scth \fa $\mu\in F_1^*$, $x\in U$, $v\in E_1^n$, 
	\[
		\mu(\vphi_n(x)(v))=[\psi^*\mu(x)]_n(v)\ .
	\]
	Thus, we have a family $\vphi=(\vphi_n)\in\mathcal O_{\mathcal U,\mathcal V}(U)$, and by \thmref{Prop}{structuresheafiso}, there exists a unique $f\in\Ct[^\infty_{MS}]0{\mathcal U,\mathcal V}$ \scth $\phi(f)=\vphi$. For $\mu\in F^*$, let $g_\mu\in\Ct[^\infty_{MS}]0{\mathcal V,\mathcal R}$ be determined by $\phi(g_\mu)=(\mu_0,\mu_1)$ where $\mu_0=\mu|_{F_0}$ and $\mu_1$ is the constant map $U\to\Hom0{F_1,\Pi R}:x\mapsto\Pi\mu|_{F_1}$. Then $\phi(g_\mu\circ\vphi)=\psi^*(\mu)$ by \thmref{Prop}{composfmla}. By applying parts (ii) and (iii) of \thmref{Lem}{genhadamard}, we find that $(\psi,\psi^*)=\Phi(f)$. 
\end{proof}

Along the way, we have also proved the following theorem.

\begin{Th}[msbklsdomequiv]
	Let $R$ be an arbitrary unital commutative Hausdorff topological $\rats$-algebra with dense group of units. Then $\SDom_{MS}$ is equivalent to the following category: objects are pairs $(U,E)$ where $E\in\TopSMod_R$ and $U\subset E_0$ is open; morphisms $(U,E)\to(V,F)$ are the elements of $\mathcal O_{\underline E_U,\underline F_V}(U)$, \emph{i.e.}~families $(\vphi_n)$ where $\vphi_n:U\to\Alt[^n]0{E_1,F_n}$ are smooth over $R$ ($F_n=F_0$ or $F_n=F_1$, according to the parity of $n$); composition is given by \eqref{eq:skelcompos}, and this determines the identity morphisms uniquely. 
\end{Th}

\begin{Rem}
	This is already stated in \cite{molotkov} (for $R=\reals$). Molotkov calls the morphisms in the category defined in \thmref{Th}{msbklsdomequiv} `skeletons'. 
\end{Rem}

\subsection{Supermanifolds as sheaves}\label{section32}

In what follows, we assume that $R$ is a non-discrete Hausdorff topological field of characteristic zero. 

We have seen the equivalence of different concepts of (finite-dimensional) superdomains. We will presently define corresponding categories of supermanifolds. Our main task will then be to show that these categories are, again, equivalent. Since we have already seen this on the level of their local models, the proof is a matter of gluing local pieces. 

More precisely, we will embed supermanifolds (in their different incarnations) into sheaves on the corresponding categories of superdomains. The equivalence of the categories of superdomains induces an equivalence of the categories of sheaves, and this restricts to an equivalence of the categories of supermanifolds. To effect this procedure, we shall need a little bit of terminology concerning sites. We refer the reader to \cite{vistoli}, \cite{giraud}. 

\subsubsection{Molotkov--Sachse supermanifolds as sheaves}

\begin{Defn}
	We define a Grothendieck topology on $\Top^\Lambda$ as follows. A natural transformation $\mathcal F'\Rightarrow\mathcal F$ is an \emph{open embedding} if $\mathcal F'(\lambda)\to\mathcal F(\lambda)$ is an open embedding of topological spaces \fa $\lambda$. Equivalently, it factors as natural equivalence of $\mathcal F'$ through an open subfunctor of $\mathcal F$. Now, we call a family $(f_\alpha:\mathcal F_\alpha\Rightarrow\mathcal F)$ of open embeddings a \emph{covering} if it is jointly surjective, \emph{i.e.}~$\bigcup_\alpha f_{\alpha\lambda}(\mathcal F_\alpha(\lambda))=\mathcal F(\lambda)$ \fa $\lambda$. 
	
	A morphism $f:\mathcal U\Rightarrow\mathcal V$ is called an \emph{open embedding of superdomains} if it factors as an isomorphism in $\SDom_{MS}$ of $\mathcal U$ with an open subfunctor of $\mathcal V$. We call a family $(f_\alpha:\mathcal U_\alpha\Rightarrow\mathcal V)$ of open embeddings of superdomains a \emph{covering} in $\SDom_{MS}$ if it is a covering in $\Top^\Lambda$. We call both topologies the \emph{DeWitt topology}. 
\end{Defn}

\begin{Rem}
	In passing, note the following subtle point: A bijective smooth map with invertible differential may not have a smooth inverse if the base field $R$ is not $\reals$ or $\cplxs$. 
\end{Rem}

\begin{Par}
	In proving that the above definition actually gives Grothendieck to\-po\-lo\-gies on $\Top^\Lambda$ and $\SDom_{MS}$, there is only one fact that is slightly non-trivial. Indeed, let $(f_\alpha:\mathcal F_\alpha\Rightarrow\mathcal F)$ be a covering and $g:\mathcal G\Rightarrow\mathcal F$ be any morphism. For any open subfunctor $\mathcal H\subset\mathcal F$, let $g^{-1}\mathcal H$ be the open subfunctor of $\mathcal G$ given by $(g^{-1}\mathcal H)(\lambda)=g_\lambda^{-1}\mathcal H(\lambda)$ for any $\lambda\in\Lambda$, and $(g^{-1}\mathcal H)(\alpha)=\mathcal G(\alpha)|_{(g^{-1}\mathcal H)(\lambda)}$ for any morphism $\alpha:\lambda\to\lambda'$. The latter is well-defined by the naturality of $g$. 
	
	Factor $f_\alpha$ into isomorphisms $f_\alpha':\mathcal F_\alpha\Rightarrow\mathcal F_\alpha'$ and inclusions of open subfunctors $\mathcal F_\alpha'\subset\mathcal F$. Define $\mathcal H_\alpha=g^{-1}\mathcal F_\alpha'$ and let $p_{\alpha1}=f_\alpha^{\prime-1}\circ g|_{\mathcal H_\alpha}$, and $p_{\alpha2}$ be the inclusion $\mathcal H_\alpha\to\mathcal G$. Then $\mathcal H_\alpha$, with the projections $p_{\alpha j}$, satisfies the universal property of the fibred product $\mathcal F_\alpha\times_{\mathcal F}\mathcal G$, and the family $(p_{\alpha2}:\mathcal H_\alpha\Rightarrow\mathcal G)$ is a covering of $\mathcal G$. 
	
	We shall always choose our fibred products in the above fashion whenever $f_\alpha$ is an open embedding (in general, we will not be able to do so). 
\end{Par}

Now we are ready to define supermanifolds modeled on Molotkov--Sachse superdomains. We will closely follow the original definitions of Molotkov--Sachse \cite{molotkov,sachse-diss}. In general, we will neglect the Hausdorff axiom, but this is easily remedied, as we shall presently see.  

\begin{Defn}
	Let $\mathcal M\in\Top^\Lambda$. A covering $\mathcal A=(\vphi_\alpha:\mathcal U_\alpha\Rightarrow\mathcal M)$ \scth the fibred products $\mathcal U_{\alpha\beta}=\mathcal U_\alpha\times_{\mathcal M}\mathcal U_\beta$ exist in $\SDom_{MS}$ is called a \emph{supermanifold atlas}. Since the object $\mathcal U_{\alpha\beta}$ lies in $\SDom_{MS}$, the latter requirement means that the projections $\mathcal U_{\alpha\beta}\Rightarrow\mathcal U_\alpha$ and $\mathcal U_{\alpha\beta}\Rightarrow\mathcal U_\beta$ are morphisms in $\SDom_{MS}$. (More precisely, for our choice of fibred products, the first of these projections is required to be a morphism in $\SDom_{MS}$.) Given another supermanifold atlas $\mathcal B$, $\mathcal A$ and $\mathcal B$ are called \emph{equivalent} if $\mathcal A\cup\mathcal B$ is a supermanifold atlas. 
	
	A pair $(\mathcal M,[\mathcal A])$ where $[\mathcal A]$ is an equivalence class of atlases on $\mathcal M$, is called a \emph{supermanifold} (in the sense of Molotkov--Sachse). We usually suppress $[\mathcal A]$ from the notation; moreover, if we say that we wish to consider an atlas of a given supermanifold, then we will always mean an atlas which belongs to the given equivalence class. A morphism of supermanifolds is a morphism $f:\mathcal M\Rightarrow\mathcal N$ in $\Top^\Lambda$ \scth for some given atlases $\mathcal A=(\vphi_\alpha:\mathcal U_\alpha\Rightarrow\mathcal M)$, $\mathcal B=(\psi_\beta:\mathcal V_\beta\Rightarrow\mathcal N)$, the pullback $\mathcal U_\alpha\times_{\mathcal N}\mathcal V_\beta$ of $f\circ\vphi_\alpha$ and $\psi_\beta$ lies with its projections in $\SDom_{MS}$. By our above considerations, this means that $\psi_\beta^{-1}\circ f\circ\vphi_\alpha:(f\circ\vphi_\alpha)^{-1}\psi_\beta(\mathcal V_\beta)\to\mathcal V_\beta$ is smooth.
\end{Defn} 
	
\begin{Rem}
	This definition of morphisms is independent of the choice of (equivalent) atlases, since smoothness of natural transformations of superdomains is a local property \cite[Lemma 4.9]{BGlN}. (This fact relies on the assumption that $R$ is a field.) 
\end{Rem}

\begin{Defn}
	We call any $\mathcal M\in\Top^\Lambda$ \emph{Hausdorff} whenever the diagonal morphism $\delta:\mathcal M\Rightarrow\mathcal M\times\mathcal M$ is closed in the sense that all its constituents $\delta_\lambda$ are. We denote the category of Molotkov--Sachse supermanifolds and their morphisms by $\SMan_{MS}=\SMan_{MS}(R)$; the full subcategory of Hausdorff supermanifolds is denoted by $\SMan_{MS}^{Hd}=\SMan_{MS}^{Hd}(R)$. Clearly, $\SDom_{MS}$ is a full subcategory of $\SMan_{MS}^{Hd}$. 
\end{Defn}

The following proposition shows that one can glue (morphisms of) superdomains to (morphisms of) supermanifolds. 

\begin{Prop}[msglue]
	Let $(\mathcal U_\alpha)\subset\SDom_{MS}$, $\mathcal U_{\alpha\beta}\subset\mathcal U_\alpha$ be open subfunctors, and $\vphi_{\beta\alpha}:\mathcal U_{\alpha\beta}\Rightarrow\mathcal U_{\beta\alpha}$ be isomorphisms in $\SDom_{MS}$ \scth we have $\vphi_{\alpha\alpha}=\id$ on $\mathcal U_{\alpha\alpha}=\mathcal U_\alpha$ and $\vphi_{\gamma\alpha}=\vphi_{\gamma\beta}\circ\vphi_{\beta\alpha}$ on $\mathcal U_{\alpha\beta}\cap\mathcal U_{\alpha\gamma}$. Then there exists a functor $\mathcal M$ with a supermanifold atlas $\vphi_\alpha:\mathcal U_\alpha\Rightarrow\mathcal M$ \scth $\vphi_\beta\circ\vphi_{\beta\alpha}=\vphi_\alpha$ on $\mathcal U_{\alpha\beta}$. Moreover, $\mathcal M$ is unique with this property, up to unique isomorphism, and $\mathcal M$ is Hausdorff if and only if $\mathcal M(R)$ is Hausdorff. 
\end{Prop}

\begin{proof}
	For each $\lambda$, define $\mathcal M(\lambda)=\coprod_\alpha\mathcal U_\alpha(\lambda)/\sim$ where the relation $\sim$ on $\coprod_{\alpha,\beta}\mathcal U_{\alpha\beta}(\lambda)$ is the union of the graphs of the isomorphisms $\vphi_{\alpha\beta,\lambda}$. By defi\-nition of the topology, there are open embeddings $\vphi_{\alpha,\lambda}:\mathcal U_{\alpha}(\lambda)\to\mathcal M(\lambda)$ which satisfy $\vphi_{\beta,\lambda}\circ\vphi_{\beta\alpha,\lambda}=\vphi_{\alpha,\lambda}$ on $\mathcal U_{\alpha\beta}(\lambda)$ (the openness follows from the fact that the saturation of $\mathcal U_\alpha(\lambda)$ is the open subset $\coprod_\beta\mathcal U_{\alpha\beta}(\lambda)$ of $\coprod_\beta\mathcal U_\beta(\lambda)$). 
	
	Using the naturality of $\vphi_{\beta\alpha}$, one readily shows that $\lambda\mapsto\mathcal M(\lambda)$ defines a functor $\mathcal M\in\Top^\Lambda$ with a supermanifold atlas as specified. (On morphisms $\phi:\lambda\to\lambda'$, $\mathcal M(\phi)$ is defined by $\mathcal M(\phi)\circ\vphi_{\alpha,\lambda}=\vphi_{\alpha,\lambda'}\circ\mathcal U_\alpha(\phi)$.) Since $\mathcal M$ is the solution of a universal problem, it is unique up to unique isomorphism. 
	
 	Finally, if $\mathcal M$ is Hausdorff, then $\mathcal M(R)$ is Hausdorff. Conversely, assume that $\mathcal M(R)$ is Hausdorff. We will use the characterisation of Hausdorff equivalence relations from \cite[\S~8.3]{bourbaki-top1}. Let $\lambda\in\Lambda$ and $x,y\in\coprod_\alpha\mathcal U_\alpha(\lambda)$ be inequivalent points. Let $x_R,y_R\in\coprod_\alpha\mathcal U_\alpha(R)$ be their respective images under $(\coprod_\alpha\mathcal U_\alpha)(\eps)$. If $x_R,y_R$ are inequivalent, then there exist by the assumption on $\mathcal M(R)$ disjoint saturated open neighbourbourhoods $U_R,V_R\subset\coprod_\alpha\mathcal U_\alpha(R)$ of $x_R$ and $y_R$, respectively. (Here, \emph{saturated} means saturated with respect to $\sim$.) By \thmref{Prop}{dewittgrothtop}, there exist unique superdomains $\mathcal U,\mathcal V$ with $\mathcal U(R)=U_R$ and $\mathcal V(R)=V_R$. It is easy to check that $\mathcal U(\lambda)$, $\mathcal V(\lambda)$ are disjoint saturated open neighbourhoods of $x$ and $y$, respectively. Thus, $\mathcal M(\lambda)$ is Hausdorff, and since $\lambda$ was arbitrary, $\mathcal M$ is Hausdorff. 
\end{proof}

\begin{Prop}[msyoneda]
	The Yoneda embedding 
	\[
	Y=Y_{MS}:\SMan_{MS}\to\Sh0{\SDom_{MS}}:\mathcal M\mapsto\mathcal M(-)=\Hom0{-,\mathcal M}
	\]
	is fully faithful. 
\end{Prop}

\begin{proof}	
	By \thmref{Prop}{msglue}, we can glue superdomains to obtain supermanifolds, and similarly, it follows easily that we can also glue morphisms of superdomains to obtain morphisms of supermanifolds. This shows that for any supermanifold $\mathcal M$, $Y(\mathcal M)$ is indeed a sheaf on $\SDom_{MS}$ with the DeWitt topology. Moreover, by the same line of thought, the presheaf $\mathcal N\mapsto\Hom0{\mathcal N,\mathcal M}$ on $\SMan_{MS}$ is entirely determined by its restriction $Y(\mathcal M)$ to $\SDom_{MS}$. Thus, the assertion follows from the usual Yoneda lemma.
\end{proof}

By the two previous propositions, we have embedded $\SMan_{MS}$ as a full subcategory of $\Sh0{\SDom_{MS}}$. To characterise the essential image of this embedding, we introduce some terminology which is common in the theory of stacks \cite{metzler}. Since we will use these concepts in a restricted setting, we will not give the most general definitions. 

\begin{Defn}
	Let $\mathcal C$ be a site and $\Presh0{\mathcal C}=\Sets^{\mathcal C^{op}}$. With any $C\in\mathcal C$, we associate the Yoneda functor $h_C=\Hom0{-,C}\in\Presh0{\mathcal C}$. Similarly, with any morphism $f:C\to C'$ in $\mathcal C$, we associate the natural transformation $h_f=\Hom0{-,f}:h_C\Rightarrow h_{C'}$. Recall that for $F\in\Presh0{\mathcal C}$, there is a natural bijection of sets $\Hom0{h_C,F}\cong F(C)$, by the token of which we will identify any $d\in F(C)$ with the associated natural transformation.
	 
	Let $f:F\Rightarrow G$ be a morphism in the category $\Sh0{\mathcal C}$ of sheaves on $\mathcal C$. We call $f$ a \emph{covering morphism} if for any $C\in\mathcal C$ and any $d\in G(C)$, there exist a covering $(f_\alpha:C_\alpha\to C)$ in $\mathcal C$ and elements $d_\alpha\in F(C_\alpha)$ \scth $f\circ d_\alpha=d\circ h_{f_\alpha}$ as natural transformations $h_{C_\alpha}\Rightarrow G$ \fa $\alpha$---equivalently, we may require that $G(f_\alpha)(d)=f_{C_\alpha}(d_\alpha)$ \fa $\alpha$. 
	
	Let $\mathcal C=\SDom_{MS}$ with the DeWitt topology, and $f:F\Rightarrow G$ be a morphism where $F,G\in\Sh0{\SDom_{MS}}$. The morphism $f$ is \emph{representable} if for any $\mathcal U\in\SDom_{MS}$ and any $g\in G(\mathcal U)$, there are $\mathcal U_\alpha\in\SDom_{MS}$ and a natural equivalence $Y\Parens1{\coprod_\alpha\mathcal U_\alpha}\cong Y(\mathcal U)\times_GF$. (Here, recall that the fibred product $F\times_HG$ in $\Sh0{\mathcal C}$ is given by the fibred product in sets.) The sheaf $F$ is \emph{locally representable} if there are $\mathcal U_\alpha\in\SDom_{MS}$ and a representable covering $Y\Parens1{\coprod_\alpha\mathcal U_\alpha}\Rightarrow F$. 
	
	We call a representable morphism $f$ an \emph{\'etale} if \fa $\mathcal U\in\SDom_{MS}$ and any representable morphism $Y(\mathcal U)\Rightarrow G$, the projection $F\times_GY(\mathcal U)\Rightarrow Y(\mathcal U)$ is (the image under $Y$ of) an \'etale morphism of supermanifolds. Here, a morphism $g:\mathcal M\Rightarrow\mathcal N$ in $\SMan_{MS}$ is called an \emph{\'etale morphism of supermanifolds} if there exists an atlas $(\vphi_\alpha:\mathcal U_\alpha\Rightarrow\mathcal M)$ \scth \fs (equivalently, any) atlas $(\psi_\beta:\mathcal V_\beta\Rightarrow\mathcal N)$ and any $\alpha,\beta$, the projection $p_2:\mathcal U_\alpha\times_{\mathcal N}\mathcal V_\beta\Rightarrow\mathcal V_\beta$ is an open embedding. Thus, there exists an open subfunctor $\mathcal U_{\alpha\beta}\subset\mathcal U_\alpha$ \scth $p_1:\mathcal U_\alpha\times_{\mathcal N}\mathcal V_\beta\Rightarrow\mathcal U_{\alpha\beta}$ is an isomorphism in $\SDom_{MS}$. If, furthermore, $\mathcal V_{\beta\alpha}\subset\mathcal V_\beta$ denotes the open subfunctor which, as a superdomain, is isomorphic to $\mathcal U_\alpha\times_{\mathcal N}\mathcal V_\beta$ via the morphism induced by $p_2$, then the local expression $g_{\beta\alpha}:\mathcal U_{\alpha\beta}\Rightarrow\mathcal V_{\beta\alpha}$ of $g$, determined by $\psi_\beta\circ g_{\beta\alpha}=g\circ\vphi_\alpha$ on $\mathcal U_{\alpha\beta}$,  is an isomorphism in $\SDom_{MS}$. 
\end{Defn}

\begin{Rem}	
	Some comments on the above terminology are perhaps in order for the unaccustomed reader. The term \emph{covering morphism} should be understood as a generalisation of the corresponding notion from topology (\emph{i.e.}~of covering maps); another analogous notion is that of a surjective submersion. 
		
	In general, `local' properties of morphisms of sheaves (such as `\'etale') can be defined only for such morphisms which are given as `glued morphisms of superdomains'. The notion `representable' encodes the concept of a morphism being in this sense `glued from local pieces'. Representable morphisms then are \'etale when so are their local representatives. 

	The reader should observe that the notion of \emph{\'etale} defined for morphisms of supermanifolds is appropriate in this context (and in this context only). It could also be referred to by the somewhat cumbersome parlance `local diffeomorphism'. 
\end{Rem}

\begin{Prop}[msrepresentable]
	A sheaf $F\in\Sh0{\SDom_{MS}}$ belongs to the essential image of $Y=Y_{MS}$ if and only if it is locally representable by a representable \'etale covering morphism $p:\coprod_\alpha Y(\mathcal U_\alpha)\Rightarrow F$ where $\mathcal U_\alpha\in\SDom_{MS}$. 
\end{Prop}

\begin{proof}
	For the sake of simplicity, we will write $Y(\mathcal U)=\mathcal U$ whenever no confusion is possible. Let $\mathcal M$ be a supermanifold and $(\vphi_\alpha:\mathcal U_\alpha\Rightarrow\mathcal M)$ be an atlas, and define $p=\amalg_\alpha \vphi_\alpha:\coprod_\alpha\mathcal U_\alpha\Rightarrow\mathcal M$. Since $\coprod_\alpha\mathcal U_\alpha$ is a supermanifold with an atlas given by the inclusions $i_\alpha:\mathcal U_\alpha\Rightarrow\coprod_\beta\mathcal U_\beta$, we may think of $p$ as a morphism in $\SMan_{MS}$ by \thmref{Prop}{msyoneda}. We need to see that $p$ is a representable \'etale covering.
	
	First, we show that it is a covering in $\Sh0{\SDom_{MS}}$. To that end, let $d:\mathcal U\Rightarrow\mathcal M$ be a morphism of supermanifolds where $\mathcal U\in\SDom_{MS}$. Since $\id:\mathcal U\Rightarrow\mathcal U$ is an atlas of $\mathcal U$, this implies that the fibre product $\mathcal U_\alpha\times_{\mathcal M}\mathcal U$ exists with its projections in $\SDom_{MS}$. Since $\vphi_\alpha$ is an open embedding in $\Top^\Lambda$, it follows that the projections $(f_\alpha:\mathcal U_\alpha\times_{\mathcal M}\mathcal U\Rightarrow\mathcal U)$ form a covering in $\SDom_{MS}$. Let $p_\alpha:\mathcal U_\alpha\times_{\mathcal M}\mathcal U\Rightarrow\mathcal U_\alpha$ be the first projections; then $d_\alpha=i_\alpha\circ p_\alpha$ satisfies $p\circ d_\alpha=\vphi_\alpha\circ p_\alpha=d\circ f_\alpha$, as required. Hence, $p$ is a covering morphism. 
	
	Next, we show that $p$ is representable. This follows in much the same way: for any morphism $g:\mathcal U\Rightarrow\mathcal M$ of supermanifolds where $\mathcal U\in\SDom_{MS}$, the pullback $\mathcal U_\alpha\times_{\mathcal M}\mathcal U$ is a superdomain. Hence, $p$ is representable. Moreover, the projections $(\mathcal U_\alpha\times_{\mathcal M}\mathcal U\Rightarrow\mathcal U)$ form a covering of superdomains, so that their disjoint union $\coprod_\alpha\mathcal U_\alpha\times_\mathcal M\mathcal U\Rightarrow\mathcal U$ is an \'etale morphism of supermanifolds, and $p$ is a representable \'etale covering morphism. 
	
	Conversely, assume that $p:\coprod_\alpha\mathcal U_\alpha\Rightarrow F$ is a representable \'etale covering where $\mathcal U_\alpha\in\SDom_{MS}$. Define $\vphi_\alpha:\mathcal U_\alpha\Rightarrow F$ by $\vphi_\alpha=p\circ i_\alpha$. Let $\mathcal N=\coprod_\alpha\mathcal U_\alpha$. By assumption, the pullback $\mathcal N\times_F\mathcal N$ of $p$ with itself exists in $\SMan_{MS}$, and the projection $g:\mathcal N\times_F\mathcal N\Rightarrow\mathcal N$ is \'etale. 
	
	Hence, there exist an atlas $(\psi_j:\mathcal V_j\Rightarrow\mathcal N\times_F\mathcal N)$ and open subfunctors $\mathcal V_j^\alpha\subset\mathcal V_\gamma$ and $\mathcal U_\alpha^j\subset\mathcal U_\alpha$ \scth the local expression $g_\alpha^j:\mathcal V_\gamma^\alpha\Rightarrow\mathcal U_\alpha^j$ of $g$, determined by $g\circ\psi_j=\vphi_\alpha\circ g_\alpha^j$ on $\mathcal V_j^\alpha$, is an isomorphism in $\SDom_{MS}$, and $\mathcal U_\alpha=\bigcup_j\mathcal U_\alpha^j$. For $\alpha,\beta,j$, let $\mathcal U_{\alpha\beta}^j\subset\mathcal U_\alpha^j$ be an open subfunctor \scth the restriction $g_\alpha^j:\mathcal V_j^\alpha\cap\mathcal V_j^\beta\Rightarrow\mathcal U_{\alpha\beta}^j$ is an isomorphism in $\SDom_{MS}$. Let $\vphi_{\beta\alpha}^j:\mathcal U_{\alpha\beta}^j\Rightarrow\mathcal U_{\beta\alpha}^j$ be defined by $\vphi_{\beta\alpha}^j\circ g_\alpha^j=g_\beta^j$. If $j'$ is another index in the same index set as $j$, then on $\mathcal U_{\alpha\beta}^j\cap\mathcal U_{\alpha\beta}^{j'}$, we have
	\[
		\vphi_\beta\circ\vphi_{\beta\alpha}^j=\vphi_\beta\circ g_\beta^j\circ(g_\alpha^j)^{-1}=g\circ\psi_j\circ(g_\alpha^j)^{-1}=\vphi_\alpha=\vphi_\beta\circ\vphi_{\beta\alpha}^{j'}\ .
	\]
	Hence, $\vphi_{\beta\alpha}^j=\vphi_{\beta\alpha}^{j'}$ on $\mathcal U_{\alpha\beta}^j\cap\mathcal U_{\alpha\beta}^{j'}$. Let $\mathcal U_{\alpha\beta}=\bigcup_j\mathcal U_{\alpha\beta}^j$; this defines an open subfunctor of $\mathcal U_\alpha$. As a special case of \thmref{Prop}{msglue}, there exists a unique isomorphism $\vphi_{\beta\alpha}:\mathcal U_{\alpha\beta}\Rightarrow\mathcal U_{\beta\alpha}$ in $\SDom_{MS}$ \scth $\vphi_{\beta\alpha}=\vphi_{\beta\alpha}^j$ on $\mathcal U_{\alpha\beta}^j$. Clearly, $\mathcal U_\alpha=\mathcal U_{\alpha\alpha}$, $\vphi_{\alpha\alpha}=\id$, $\vphi_{\gamma\beta}\circ\vphi_{\beta\alpha}=\vphi_{\gamma\alpha}$ on $\mathcal U_{\alpha\beta}\cap\mathcal U_{\alpha\gamma}$. By \thmref{Prop}{msglue}, there exists $\mathcal M\in\SMan_{MS}$ and an atlas $(\tilde\vphi_\alpha:\mathcal U_\alpha\Rightarrow\mathcal M)$ \scth $\tilde\vphi_\alpha\circ\vphi_{\alpha\beta}=\tilde\vphi_\beta$ on $\mathcal U_{\beta\alpha}$. 
	
	We wish to define a morphism $\vphi:\mathcal M\Rightarrow F$ in $\Sh0{\SDom_{MS}}$ \scth $\vphi_\alpha=\vphi\circ\tilde\vphi_\alpha$. Let $\mathcal V\in\SDom_{MS}$. We have to define $\vphi_{\mathcal V}:\mathcal M(\mathcal V)\to F(\mathcal V)$. To that end, let $\phi\in\mathcal M(\mathcal V)$. Set $\mathcal V_\alpha=\mathcal U_\alpha\times_{\mathcal M}\mathcal V$, let $\phi_\alpha:\mathcal V_\alpha\Rightarrow\mathcal U_\alpha$ be the first projection, and $\psi_\alpha:\mathcal V_\alpha\Rightarrow\mathcal V$ the second projection. Then $\tilde\vphi_\alpha\circ\phi_\alpha=\phi\circ\psi_\alpha$. Set $\tilde\phi_\alpha=\vphi_\alpha\circ\phi_\alpha:\mathcal V_\alpha\Rightarrow F$. Then $\tilde\phi_\alpha$ may be considered as an element of $F(\mathcal V_\alpha)$. Moreover, since $(\tilde\vphi_\alpha:\mathcal U_\alpha\Rightarrow\mathcal M)$ is a covering in $\SMan_{MS}$, $(\psi_\alpha:\mathcal V_\alpha\Rightarrow\mathcal V)$ is a covering in $\SDom_{MS}$. By construction, $\mathcal U_{\alpha\beta}$ is the fibre product $\mathcal U_\alpha\times_{\mathcal M}\mathcal U_\beta$ in $\SMan_{MS}$ with first projection the inclusion $i_{\alpha\beta}:\mathcal U_{\alpha\beta}\subset\mathcal U_\alpha$, and second projection $i_{\beta\alpha}\circ\vphi_{\beta\alpha}:\mathcal U_{\alpha\beta}\Rightarrow\mathcal U_\beta$. 
	
	Consider the fibre product $\mathcal V_\alpha\times_{\mathcal M}\mathcal V_\beta$ with its projections $p_\alpha$, $p_\beta$ with codomains $\mathcal V_\alpha$ and $\mathcal V_\beta$, respectively. We have
	\[
		\tilde\vphi_\alpha\circ\phi_\alpha\circ p_\alpha=\phi\circ\psi_\alpha\circ p_\alpha=\phi\circ\psi_\beta\circ p_\beta=\tilde\vphi_\beta\circ\phi_\beta\circ p_\beta\ .
	\]
	By the universal property of the fibred product, there exists a unique morphism $\phi_{\alpha\beta}:\mathcal V_\alpha\times_{\mathcal V}\mathcal V_\beta\to\mathcal U_{\alpha\beta}$ \scth 
	\[
	\phi_\alpha\circ p_\alpha=i_{\alpha\beta}\circ\phi_{\alpha\beta}\nd\phi_\beta\circ p_\beta=i_{\beta\alpha}\circ\vphi_{\beta\alpha}\circ\phi_{\alpha\beta}\ .
	\]
	Then 
	\[
		\tilde\phi_\beta\circ p_\beta=\vphi_\beta\circ\phi_\beta\circ p_\beta=\vphi_\beta\circ\vphi_{\beta\alpha}\circ\phi_{\alpha\beta}=\vphi_\alpha\circ\phi_{\alpha\beta}=\tilde\phi_\alpha\circ p_\alpha\ .
	\]
	Since $F$ is a sheaf, there exists a unique $\tilde\phi\in F(\mathcal V)$ \scth $\tilde\phi\circ\psi_\alpha=\tilde\phi_\alpha$ \fa $\alpha$. We set $\vphi_{\mathcal V}(\phi)=\tilde\phi$. Summarising the definition, $\phi$ and $\tilde\phi$ are related by $\tilde\phi\circ\psi_\alpha=\vphi_\alpha\circ\phi_\alpha$ and $\tilde\vphi_\alpha\circ\phi_\alpha=\phi\circ\psi_\alpha$. One readily checks that this defines a natural transformation $\vphi:\mathcal M\Rightarrow F$ in $\Sh0{\SDom_{MS}}$ and that $\vphi\circ\tilde\vphi_\alpha=\vphi_\alpha$. The claim is proved as soon as it has been shown that $\vphi$ is an isomorphism.
	
	To that end, it suffices to prove that the natural inclusion of fibred products $\coprod_\alpha\mathcal U_\alpha\times_\mathcal M\coprod_\alpha\mathcal U_\alpha\to\coprod_\alpha\mathcal U_\alpha\times_F\coprod_\alpha\mathcal U_\alpha$ is an isomorphism. But this follows from the fact that $\vphi$ is an isomorphism on the open subfunctor of $\mathcal M$ onto which $\tilde\vphi_\alpha$ defines an isomorphism. 
\end{proof}

\subsubsection{Berezin--Kostant--Leites supermanifolds as sheaves}

We will now apply the same procedure that we have applied to Molotkov--Sachse supermanifolds in the context of Berezin--Kostant--Leites supermani\-folds over $R$. \emph{I.e.}, we define Berezin--Kostant--Leites supermanifolds, embed them into the category of sheaves on Berezin--Kostant--Leites superdomains, and characterise the essential image of this embedding. Most of the arguments will be completely parallel to the case of Molotkov--Sachse supermanifolds, so that we will only indicate the differences. 

\begin{Defn}
	Let $\SRSp_R$ denote the category whose objects (called \emph{graded $R$-ringed spaces}) are pairs $(X,\mathcal F)$ where $X$ is a topological space, and $\mathcal F$ is a sheaf of unital $R$-superal\-ge\-bras over $X$; and whose morphisms $(X,\mathcal F)\to(Y,\mathcal G)$ are pairs $(\vphi,\vphi^*)$ where $\vphi:X\to Y$ is continuous and $\vphi^*:\mathcal G\to\vphi_*\mathcal F$ is a morphism of sheaves of unital $R$-superalgebras. Clearly, $\SDom_{BKL}$ is a full subcategory. 
	
	A morphism $(\vphi,\vphi^*):(X,\mathcal F)\to(Y,\mathcal G)$ is called an \emph{open embedding} if $\vphi$ is an open embedding and $\vphi^*$ induces an isomorphism $\mathcal G|_{\vphi(X)}\to\vphi_*\mathcal F$. A family of open embeddings in $\SRSp_R$ is a \emph{covering} if the underlying family of maps of topological spaces is  jointly surjective. Clearly, this defines a Grothendieck topology on $\SDom_{BKL}$ and $\SRSp_R$ which we call the \emph{standard topology}. 
	
	An object $(X,\mathcal F)\in\SRSp_R$ is a \emph{supermanifold} in the sense of Berezin--Kostant--Leites if it has a covering by objects of $\SDom_{BKL}$. We define $\SMan_{BKL}$ to be the full subcategory of $\SRSp_R$ whose objects are supermanifolds. 
	
	We call an object $\mathcal X=(X,\mathcal F)\in\SRSp_R$ \emph{Hausdorff} if $X$ is Hausdorff; equivalently, the diagonal morphism $\mathcal X\to\mathcal X\times\mathcal X$ is closed (on the level of the underlying topological spaces). We let $\SMan_{BKL}^{Hd}$ be the full subcategory of $\SMan_{BKL}$ whose objects are Hausdorff. If $R=\reals$ or $R=\cplxs$, then $\SMan_{BKL}^{Hd}$ is the category of Berezin--Kostant--Leites supermanifolds as it is usually defined in the literature (if, with a view towards infinite dimensions, one ignores the axiom of second countability). 
\end{Defn}

\begin{Prop}[bklyoneda]
	The Yoneda embedding
	\[
		Y=Y_{BKL}:\SMan_{BKL}\to\Sh0{\SDom_{BKL}}:\mathcal X\mapsto\mathcal X(-)=\Hom0{-,\mathcal X}
	\]
	is fully faithful. 
\end{Prop}

\begin{proof}
	This follows from the Yoneda lemma and the fact that ringed spaces and their morphisms can be glued \cite[Chapitre 0, (4.1.7)]{grothendieck-ega1}.
\end{proof}

We define the concepts of representable and \'etale morphisms in the category $\Sh0{\SDom_{BKL}}$ in much the same way as for $\Sh0{\SDom_{MS}}$. 

\begin{Defn}
	Consider $\SDom_{BKL}$ with the standard topology, and let $f:F\Rightarrow G$ be a morphism where $F,G\in\Sh0{\SDom_{BKL}}$. The morphism $f$ is \emph{representable} if for any $\mathcal U\in\SDom_{BKL}$ and any $g\in G(\mathcal U)$, there are superdomains $\mathcal U_\alpha\in\SDom_{BKL}$ and a natural equivalence $Y\Parens1{\coprod_\alpha\mathcal U_\alpha}\cong Y(\mathcal U)\times_GF$. The sheaf $F$ is \emph{locally representable} if there are $\mathcal U_\alpha\in\SDom_{BKL}$ and a representable covering $Y\Parens1{\mathcal U_\alpha}\Rightarrow F$. 
	
	We call $f$ an \emph{\'etale morphism} if \fa $\mathcal U\in\SDom_{BKL}$ and any representable morphism $Y(\mathcal U)\Rightarrow G$, the projection $F\times_GY(\mathcal U)\Rightarrow Y(\mathcal U)$ is (the image under $Y$ of) an \'etale morphism of supermanifolds. Here, a morphism $g:\mathcal M\to\mathcal N$ in $\SMan_{BKL}$ is called an \emph{\'etale morphism of supermanifolds} if there exists a covering $(\vphi_\alpha:\mathcal U_\alpha\to\mathcal M)$ of $\mathcal M$ by objects of $\SDom_{BKL}$ \scth $(g\circ\vphi_\alpha:\mathcal U_\alpha\to\mathcal N)$ is a covering. 
\end{Defn}

\begin{Rem}
	Again, the above notion of \emph{\'etale} is modeled on local homeomorphisms, which is entirely appropriate in this context. It should not be confused with the more subtle notion of \emph{\'etale} from algebraic geometry. 
\end{Rem}

\begin{Prop}[bklrepresentable]
	A sheaf $F\in\Sh0{\SDom_{BKL}}$ belongs to the essential image of $Y=Y_{BKL}$ if and only if it is locally representable by a representable \'etale covering morphism $p:\coprod_\alpha Y(\mathcal U_\alpha)\Rightarrow F$ where $\mathcal U_\alpha\in\SDom_{BKL}$. 
\end{Prop}

\begin{proof}
	The proof is exactly the same as for the setting of Molotokov--Sachse supermanifolds (\thmref{Prop}{msrepresentable}). In fact, the subtle points in that case, concerning different Grothendieck topologies, are trivial here: $\SDom_{BKL}$, $\SMan_{BKL}$ are full subcategories of $\SRSp_R$, and a covering in $\SRSp_R$ whose domains and codomains lie in $\SMan_{BKL}$ (resp.~$\SDom_{BKL})$ is a covering in that site. 
\end{proof}

\subsection{Equivalence of $\SMan_{BKL}$ and $\SMan_{MS}^{fd}$}\label{section33}

\begin{Defn}
	A supermanifold $(\mathcal M,[\mathcal A])\in\SMan_{MS}$ is called \emph{finite-dimensional} if there exists a supermanifold atlas $(\vphi_\alpha:\mathcal U_\alpha\Rightarrow\mathcal M)\in[\mathcal A]$ \scth $\mathcal U_\alpha\in\SDom_{MS}^{fd}$ \fa $\alpha$. The category $\SMan_{MS}^{fd}$ (resp.~$\SMan_{MS}^{fd,Hd}$) is defined to be the full subcategory of $\SMan_{MS}$ (resp.~$\SMan_{MS}^{Hd}$) whose objects are finite-dimensional supermanifolds. 
\end{Defn}

\begin{Rem}
	In passing, note the following amusing fact. In the definition of finite-dimensional Molotkov--Sachse supermanifolds, we have not imposed a bound on the dimension of their local models. Indeed, there is no reason to do so. For instance, if one defines (finite-dimensional) smooth manifolds properly (and not only manifolds of pure dimension), there may be an infinite number of connected components, and also an unbounded family of local dimensions. Thus, $\coprod_{n\in\nats}\reals^n$ is an ordinary smooth manifold (its global dimension is $\infty$, but all local dimensions are finite). 
\end{Rem}

\begin{Par}	
	By Propositions \ref{prop:msyoneda} and \ref{prop:msrepresentable} one clearly has by restriction a Yoneda embedding $Y=Y_{MS}^{fd}:\SMan_{MS}^{fd}\to\Sh0{\SDom_{MS}^{fd}}$ whose essential image consists of those $F\in\Sh0{\SDom_{MS}^{fd}}$ which are locally representable by a representable \'etale covering $p:\coprod_\alpha\mathcal U_\alpha\Rightarrow F$ where $\mathcal U_\alpha\in\SDom_{MS}^{fd}$ \fa $\alpha$. 
\end{Par}

\begin{Th}[msbklfdequiv]
	Let $R$ be a non-discrete Hausdorff topological field of characteristic zero, and recall the equivalence $\Phi:\SDom_{MS}^{fd}\to\SDom_{BKL}$ from \thmref{Th}{sdomequiv}. Define 
	\[
		\Phi^*:\Sh0{\SDom_{BKL}}\to\Sh0{\SDom_{MS}^{fd}}:F\mapsto F\circ\Phi\ .
	\]
	Then $\Phi^*$ is an equivalence of categories which restricts to equivalences 
	\[
	\SMan_{BKL}\to\SMan_{MS}^{fd}\nd\SMan_{BKL}^{Hd}\to\SMan_{MS}^{fd,Hd}\ .
	\]
\end{Th}

\begin{proof}
	To see that $\Phi^*$ is well-defined, we need to check that for any family of morphisms $(f_\alpha:\mathcal U_\alpha\Rightarrow\mathcal V)$ in $\SDom_{MS}^{fd}$, $(\Phi(f_\alpha):\Phi(\mathcal U_\alpha)\to\Phi(\mathcal V))$ is a covering in $\SDom_{BKL}$ if and only if $(f_\alpha)$ is a covering in $\SDom_{MS}^{fd}$. First we note that in both categories of superdomains, an open embedding was by definition the composition of an isomorphism with an embedding of subdomains. These are clearly mapped to each other by the equivalence $\Phi$. If $(f_\alpha)$ is jointly surjective, then so is $(f_{\alpha,R})$, which is the family of maps underlying the family of morphisms $(\Phi(f_\alpha))$. 
	
	Conversely, let $(\Phi(f_\alpha))$ be a covering so that $(f_\alpha)$ consists of open embeddings and we need to check that $(f_{\alpha,\lambda})$ is jointly surjective for any $\lambda\in\Lambda$. Let $x\in\mathcal V(\lambda)$, with $x_R=\mathcal V(\eps)(x)$. Then $x_R\in\Phi(f_\alpha)(\mathcal U_\alpha(R))$ \fs $\alpha$, so that by \thmref{Prop}{lambdapoints}, $x$ defines a morphism $*_\lambda\to\Phi(\mathcal V)|_{\Phi(f_\alpha)(\mathcal U_\alpha(R))}\cong\Phi(\mathcal U_\alpha)$. In other words, there is a morphism $y:*_\lambda\to\Phi(\mathcal U_\alpha)$ \scth $\Phi(f_\alpha)\circ y=x$. The corresponding point $y\in\mathcal U_\alpha(\lambda)$ satisfies $f_{\alpha,\lambda}(y)=x$. Thus, $\Phi^*$ is well-defined, and it is easy to see that it is a functor. 

	The notion of covering in a Grothendieck topology is by definition stable under isomorphism. Therefore, any quasi-inverse $\Psi$ of $\Phi$ defines in the same way a quasi-inverse $\Psi^*$ of $\Phi^*$, and $\Phi^*$ is an equivalence. 
	
	By these considerations, it is also clear that covering, representable and \'etale morphisms in one category of sheaves are mapped to the same variety of morphisms in the other category by $\Phi^*$ and $\Psi^*$, so that we obtain an equivalence $\SMan_{BKL}\to\SMan_{MS}^{fd}$. The equivalence of the categories of Hausdorff supermanifolds follows from \thmref{Prop}{msglue}, since up to isomorphism, any Molotkov--Sachse supermanifold is given as the gluing of local data as in that proposition.
\end{proof}

\section{DeWitt--Tuynman supermanifolds}\label{section4}

\subsection{Tuynman-smooth maps}

Let $R$ be a unital commutative Hausdorff topological $\rats$-algebra with a dense group of units \emph{which we assume to be locally $k_\omega$} (so that $R$ is locally compact if it is metrisable \cite[Proposition 4.8]{ggh-kac}). We denote by $\TopSMod_R^{lko}$ the full subcategory of $\TopSMod_R$ formed by the graded topological $R$-modules which are locally $k_\omega$-spaces.

The following slightly strange definition is the correct formulation of Tuynman's $\mathcal C^1$ concept \cite{tuynman-ck,tuynman} over $R$. 

\begin{Defn}
	Let $E,F\in\TopSMod_R^{lko}$, $E_\infty=\underline E(\lambda^\infty)$ and $E^\infty=E\otimes\lambda^\infty$. A subset $U\subset E^\infty$ is called \emph{nilsaturated} if $x+n\in U$ \fa $x\in U$, and all $n\in N_U=\Span0U_R\cap(E\otimes\lambda^{\infty+})$ where $\Span0U_R$ denotes the $R$-linear span. If $U$ is nilsaturated, then $U=U_R+N_U$ \fs unique $U_R\subset E$. A map $f:U\to F^\infty$ is called \emph{grounded} if $f(U_R)\subset F$; it is called \emph{even} if $U\subset E_\infty$ and $f(U)\subset F_\infty$. 
	
	If $U\subset E^\infty$ is nilsaturated, then 
	\[
	U^{[1]}=\Set1{(x,v,t)\in U\times E^\infty\times\lambda^\infty}{x+vt\in U}
	\]
	is nilsaturated in $E^{\infty[1]}=E^\infty\times E^\infty\times\lambda^\infty$.	Let $U\subset E^\infty$ be nilsaturated and $f:U\to F^\infty$ be grounded and DeWitt-continuous. (Recall the definition of the DeWitt topology from \thmref{Par}{dewitt-top}.) Then $f$ is called \emph{Tuynman-$\mathcal C^1$} whenever there exists a grounded and DeWitt-continuous map $f^{[1]}:U^{[1]}\to F^\infty$ \scth
	\[
		f(x+vt)-f(x)=f^{[1]}(x,v,t)\cdot t\mathfa(x,v,t)\in U^{[1]}\ .
	\]
	
	By induction, one defines Tuynman-$\mathcal C^k$ and Tuynman-smooth maps. 
\end{Defn}

\begin{Prop}[tuynman-taylor1]
	Let $U\subset E_\infty$ be DeWitt-open and $f:U\to F^\infty$ be even and Tuynman-$\mathcal C^2$. Then $f_R=f|_{U_R}:U_R\to F$ is $\mathcal C^1$ and there exists a DeWitt-continuous extension $f^{(1)}:U\times E^\infty\to F^\infty$ of $df_R$ \scth \fa $x\in U$, $f^{(1)}(x)=f^{(1)}(x,\cdot)$ is even and $\lambda^\infty$-linear and 
	\[
		f(x+a)-f(x)=f^{(1)}(x)(a)\mathfa x\in U\,,\,a\in\underline E(\lambda^\infty\theta_p)\ .
	\]
\end{Prop}

\begin{proof}
	We compute for $a=b\theta_p\in\underline E(\lambda^\infty\theta_p)$
	\begin{multline*}
		f(x+a)-f(x)=f^{[1]}(x,b,\theta_p)\cdot\theta_p\\
		=f^{[1]}(x,b,0)\cdot\theta_p+f^{[2]}(x,b,0,0,0,1,\theta_p)\cdot\theta_p^2=f^{(1)}(x)(b\cdot\theta_p)
	\end{multline*}
	where we set $f^{(1)}(x)(v)=f^{[1]}(x,v,0)$. By the equation, $f^{[1]}(x,v,0)$ is uniquely determined, and it is hence easy to check that $f^{(1)}$ is DeWitt-continuous, and even and $\lambda^\infty$-linear in its second argument. Since $f$, $f^{[1]}$ are assumed to be DeWitt-continuous, it follows that $f_R$, $f^{[1]}_R$ are continuous, and one sees that $f$ is $\mathcal C^1$, and that $f^{(1)}$ extends $df_R$. 
\end{proof}

This suggests a simpler definition of Tuynman-smoothness. 

\begin{Defn}
	Let $E,F\in\TopSMod_R^{lko}$ and $U\subset E_\infty$ a DeWitt-open subset. By \thmref{Lem}{dewitt-top}, $U=U_R+E_\infty^+$ \fs uniquely determined open subset $U_R\subset E_0$; here, $E_\infty^+=(E\otimes\lambda^{\infty+})_0$. 
	
	We will call any grounded and even map $f=f^{(0)}:U\to F_\infty$ \emph{$\lambda^\infty$-smooth} if $f$ is DeWitt-continuous, $f_R=f|_{U_R}:U_R\to F_0$ is $\mathcal C^\infty$ over $R$, and there exist \fa $k\sge1$ DeWitt-continuous maps $f^{(k)}:U\times(E\otimes\lambda^\infty)^k\to F\otimes\lambda^\infty$ subject to the following conditions: 
	\begin{enumerate}
		\item For $k\sge 1$ and $x\in U$, $f^{(k)}(x)$ is even, $\lambda^\infty$-multilinear, supersymmetric, and $f^{(k)}(x)(E^k)\subset F$. That $f^{(k)}(x)$ is \emph{supersymmetric} means 
		\[
			f^{(k)}(x)(\dotsc,u,v,\dotsc)=(-1)^{\Abs0u\Abs0v}f^{(k)}(x)(\dotsc,v,u,\dotsc)
		\]
		for all homogeneous $u$, $v$ of parity $\Abs0u$ and $\Abs0v$, respectively. 
		\item For any $k\sge1$, $f^{(k)}$ extends $d^kf_R$. 
		\item For any $k\sge0$, $p\sge1$, $x\in U$, $a\in\underline E(\lambda^\infty\theta_p)$, $v_1,\dotsc,v_k\in E_\infty$, we have
		\begin{equation}\label{eq:tuynmansmoothrec}
			f^{(k+1)}(x)(a,v_1,\dotsc,v_k)=\Parens1{f^{(k)}(x+a)-f^{(k)}(x)}(v_1,\dotsc,v_k)\ .
		\end{equation}
	\end{enumerate}

	It is clear that for $k\sge1$, $f^{(k)}$ is uniquely determined by its restriction to $U\times E^k\subset U_\infty\times E_\infty^k$. In particular, \eqref{eq:tuynmansmoothrec} holds \fa $v_j\in E\otimes\lambda^\infty$. 
\end{Defn}

\begin{Lem}[dersmooth]
	Let $f=f^{(0)}:U_\infty\to F_\infty$ be $\lambda^\infty$-smooth. Then \fa $k\sge0$, $f^{(k)}:U\times E_\infty^k\to F_\infty$ is continuous for the standard topology. 
\end{Lem}

\begin{proof}
	By \thmref{Prop}{dewittgrothtop}, there exists a unique open subfunctor $\mathcal U\subset\underline E$ \scth $\mathcal U(\lambda^\infty)=U$. Moreover, $U\times E_\infty^k=\varinjlim_N\mathcal U(\lambda^N)\times E_\infty^k$ in $\Top$ \cite[Proposition 4.7]{ggh-kac}. Thus, it suffices to prove the assertion for the restriction of $f^{(k)}$ to $\mathcal U(\lambda^N)\times E_\infty^k$ for any $N$. 
	
	Fix $N$. Any $x\in\mathcal U(\lambda^N)$ can be uniquely decomposed as $x=\sum_{j=0}^Nx_j$ where $x_0\in U_R$, $x_j\in\mathcal U(\theta_jR[\theta_{j+1},\dotsc,\theta_N])$. By \eqref{eq:tuynmansmoothrec}, 
	\[
		f^{(k)}(x)=f^{(k)}(x-x_N)+f^{(k+1)}(x-x_N)(x_N,\cdot)\ .
	\]
	Thus, by induction, \fa $k\sge0$, $f^{(k)}$ is continuous when restricted to $\mathcal U(\lambda^N)\times E_\infty^k$, if \fa $k\sge0$, $f^{(k)}$ is continuous when restricted to $U_R\times E_\infty^k$. 
	
	Let us prove this last assertion. As above, $U_R\times E_\infty^k=\varinjlim_MU_R\times\underline E(\lambda^M)^k$ in $\Top$, so it suffices to prove it for the restriction of $f^{(k)}$ to $U_R\times\underline E(\lambda^M)^k$ for any $M$. For $v_1,\dotsc,v_k\in\underline E(\lambda^M)$, write $v_i=\sum_Iv_{iI}\theta_I$ where $v_{iI}\in E_{\Abs0I}$, the symbol $\Abs0I$ denoting the parity of $I$. Then, \fa $x\in U_R$, 
	\[
		f^{(k)}(x)(v_1,\dotsc,v_k)=\sum_{I_1,\dotsc,I_k}f^{(k)}(x)(v_{1I_1},\dotsc,v_{kI_k})\theta_{I_1}\dotsm\theta_{I_k}\ .
	\]
	Since this sum has a fixed finite length, $f^{(k)}$ is continuous when restricted to $U_R\times\underline E(\lambda^M)^k$ if it is continuous when restricted to $U_R\times E^k$. But this is true by the definition of the DeWitt topology. 
\end{proof}

\begin{Prop}
	Let $f:U\to F_\infty$ be $\lambda^\infty$-smooth. Then \fa $x\in U$, $y\in E_\infty^+$, we have
	\begin{equation}\label{eq:tuynman-taylor}
		f(x+y)=\sum_{k=0}^\infty\frac1{k!}\cdot f^{(k)}(x)(y,\dotsc,y)\ .
	\end{equation}
\end{Prop}

\begin{proof}
	Write $y=\sum_{j=1}^Ny_j$ where $y_j\in\underline E(\lambda^\infty\theta_j)$ (such a decomposition may not be unique, but we don't care). One proves the formula 
	\[
		f(x+y)=\sum_{k=0}^\infty\sum_{\Abs0I=k}f^{(k)}(x)(y_{i_1},\dotsc,y_{i_k})
	\]
	for arbitrary $x\in U$, by induction on $N$. Then the assertion follows as in the proof of \thmref{Prop}{taylor-ms}. 
\end{proof}

\begin{Th}[tsmoothequiv]
	Let $R$ be a locally $k_\omega$ topological $\rats$-algebra with dense group of units. Let $E,F\in\TopSMod_R^{lko}$ and $\mathcal U\subset\underline E$ an open subfunctor. Let $U=\mathcal U(\lambda^\infty)$ and $f:U\to F_\infty=\underline F(\lambda^\infty)$ be a map. The following are equivalent:
	\begin{enumerate}
		\item The map $f$ is $\lambda^\infty$-smooth.
		\item The restrictions $f_\lambda=f|_{\mathcal U(\lambda)}$, $\lambda\in\Lambda$, define a smooth natural transformation $(f_\lambda):\mathcal U\Rightarrow\underline F$. 
		\item The map $f$ is Tuynman-smooth. 
	\end{enumerate}
\end{Th}

\begin{proof}
	(i) $\Rightarrow$ (ii). Let $f$ be $\lambda^\infty$-smooth. By \eqref{eq:tuynman-taylor}, and the $\lambda^\infty$-linearity and evenness of the $f^{(k)}$, $f_\lambda(\mathcal U(\lambda))\subset\underline F(\lambda)$ and it is easy to check that this defines a natural transformation of set-valued functors. By \thmref{Lem}{dersmooth}, the $f_\lambda$ are continuous. Since continuous multilinear maps are smooth (over $R$), and all of the $d^kf_R$ are smooth, \eqref{eq:tuynman-taylor} also proves that all of the $f_\lambda$ are smooth over $R$. (On $\mathcal U(\lambda)$, $\lambda\in\Lambda$ fixed, the sum in \eqref{eq:tuynman-taylor} has bounded length.) Since $df_\lambda(x)=f^{(1)}(x,\cdot)$ on their common domain of definition, $f_\lambda$ has $\lambda_0^\infty$-linear first derivatives at every point, and from \thmref{Th}{smoothequiv}, we find that $(f_\lambda):\mathcal U\Rightarrow\underline F$ is a smooth natural transformation. 
	
	(ii) $\Rightarrow$ (iii). This follows from \thmref{Prop}{taylor-ms}. 
	
	(iii) $\Rightarrow$ (i). This follows by induction from \thmref{Prop}{tuynman-taylor1}. 
\end{proof}

\begin{Prop}[tuynmantuynman]
	Let $R=\reals$ and $E,F\in\TopSMod_\reals$ be finite-dimensional. Further, let $U\subset E_\infty=\underline E(\lambda^\infty)$ be DeWitt open, and $f:U\to F_\infty$ be even, grounded and DeWitt-continuous. If $f$ is $\lambda^\infty$-smooth, then there exists a grounded DeWitt-continuous map $\phi:U\times U\times E^\infty\to F^\infty=F\otimes\lambda^\infty$ \scth $\phi(x,y)=\phi(x,y,\cdot)$ is $\lambda^\infty$-linear and even, and 
	\[
		f(x)-f(y)=\phi(x,y)(x-y)\mathfa x,y\in U\ .
	\]
	Conversely, if such a $\phi$ exists, then $f$ is Tuynman-$\mathcal C^1$. 
\end{Prop}

\begin{proof}
	Let $f$ be $\lambda^\infty$-smooth. For any $\lambda\in\Lambda^\infty$, $f_\lambda=f|_{\mathcal U(\lambda)}:\mathcal U(\lambda)\to\underline F(\lambda^\infty)$ is $\mathcal C^1$ over $\reals$, by \thmref{Th}{tsmoothequiv}. By \cite[Proposition 1.8, Remark 1.12]{tuynman}, there exist for any $\lambda\in\Lambda$ continuous maps $\phi_\lambda:\mathcal U(\lambda)\times\mathcal U(\lambda)\times\underline E(\lambda)\to\underline F(\lambda)$ \scth \[
	f_\lambda(x)-f_\lambda(y)=\phi_\lambda(x,y)(x-y)\mathfa x,y\in\mathcal U(\lambda)
\]
and $\phi_\lambda$ is $\reals$-linear in its third argument. Using the naturality of $(f_\lambda)$ and the $\lambda_0$-linearity of the derivatives, it is not hard to show that one may choose $\phi_\lambda(x,y)$ \scth an even and $\lambda$-linear extension to $E\otimes\lambda$ exists. Now, one takes $\phi=\varinjlim_\lambda\phi_\lambda$. 

	As for the converse, one may define $f^{[1]}(x,v,t)=\phi(x+vt,x)(v)$. 
\end{proof}

\begin{Rem}
	By \thmref{Prop}{tuynmantuynman}, our definition of Tuynman-smoothness is equivalent to Tuynman's definition of smoothness \cite[Definitions 1.16]{tuynman} if $R=\reals$ and the domain of definition is given by a finite-dimensional super-vector space. (Tuynman actually also considers maps which are not even, but this is clearly not a restriction, since one may take $F\oplus\Pi F$ as the target space.)
	
	Thus, our definition is an extension to arbitrary base rings $R$ (subject to our assumptions) of Tuynman's concept. A severe restriction which has to be imposed is that the model spaces $E$ are locally $k_\omega$. In case $R$ is a non-discrete locally compact field of characteristic zero, this reduces essentially to the case where $E$ is the topological direct limit of finite-dimensional subspaces. This indicates that the intuitive point of view offered by DeWitt--Tuynman's approach is not available in more general infinite-dimensional settings (such as the Banach space setting considered by Molotkov--Sachse). 
\end{Rem}

\subsection{De Witt--Tuynman vs.~Molotkov--Sachse supermanifolds}

Let $R$ be a non-discrete locally $k_\omega$ topological field of characteristic zero. 

\begin{Defn}
	The category $\SDom_{dWT}=\SDom_{dWT}(R)$ of DeWitt--Tuyn\-man superdomains has as objects pairs $(U,E)$ where $E\in\TopSMod_R^{lko}$ and $U\subset\underline E(\lambda^\infty)$ is DeWitt open; its morphisms $f:(U,E)\to(V,F)$ are the grounded, even and $\lambda^\infty$- (or, equivalently, Tuynman-) smooth mappings $f:U\to\underline F(\lambda^\infty)$ \scth $f(U)\subset V$. 
	
	Let $\SDom_{MS}^{lko}$ denote the full subcategory of $\SDom_{MS}$ of superdomains whose model spaces $E$ are locally $k_\omega$. 
\end{Defn}

\begin{Th}[dewitttuynman-ms-dom]
	The categories $\SDom_{dWT}$ and $\SDom_{MS}^{lko}$ are isomorphic. 
\end{Th}

\begin{proof}
	This follows immediately from \thmref{Th}{tsmoothequiv}. 
\end{proof}

\begin{Defn}
	Let $M$ be a topological space. Consider a jointly surjective collection $\ger A=(\vphi_\alpha:U_\alpha\to M)$ of open embeddings (for the DeWitt topology) of superdomains $U_\alpha=\mathcal U_\alpha(\lambda^\infty)$ where $\mathcal U_\alpha\subset\underline{E_\alpha}$ (for graded $R$-modules $E_\alpha\in\TopSMod_R^{lko}$) are open subfunctors. If \fa $\alpha,\beta$, 
	\[
	\vphi_{\beta\alpha}=\vphi_\beta^{-1}\circ\vphi_\alpha:(U_{\alpha\beta},E_\alpha)\to(U_{\beta\alpha},E_\beta)
	\]
	is an isomorphism in $\SDom_{dWT}$ where $U_{\alpha\beta}=\vphi_\alpha^{-1}(\vphi_\beta(U_\beta))$, then we call $\ger A$ a \emph{DeWitt--Tuynman atlas}. 
	
	If $\ger B$ is another atlas, then $\ger A$ and $\ger B$ are called \emph{equivalent} if their union is again a DeWitt--Tuynman atlas. A pair $(M,[\ger A])$ where $M$ is a topological space and $[\ger A]$ is an equivalence class of DeWitt--Tuynman atlases is called a \emph{DeWitt--Tuynman supermanifold}. We will usually not expressly mention the chosen equivalence class of atlases. 
	
	Let $f:M\to N$ be a map where $M$ and $N$ are DeWitt--Tuynman supermanifolds. Then $f$ is a morphism of DeWitt--Tuynman supermanifolds if it is continuous, and for some (any) atlases $(\vphi_\alpha)$ of $M$ and $(\psi_\beta)$ of $N$ (in the chosen equivalence classes), the map $\psi_\beta^{-1}\circ f\circ\vphi_\alpha$ is (on its domain of definition) a morphism in $\SDom_{dWT}$. 

	The topology induced by an atlas (\emph{i.e.}~the finest topology \scth all local charts $\vphi_\alpha$ are continuous for the standard topology on $U_\alpha$) will be called the \emph{standard topology} on a DeWitt--Tuynman supermanifold. It depends only on the equivalence class of the atlas, by \thmref{Lem}{dersmooth}. A DeWitt--Tuynman supermanifold will be called \emph{Hausdorff} if it is Hausdorff in the standard topology. 
	
	We denote the category of DeWitt--Tuynman supermanifolds and their morphisms by $\SMan_{dWT}=\SMan_{dWT}(R)$; the full subcategory of Hausdorff supermanifolds is denoted $\SMan_{dWT}^{Hd}$. Let $\SMan_{MS}^{lko}$ and $\SMan_{MS}^{lko,Hd}$ denote the full subcategories of $\SMan_{MS}$ and $\SMan_{MS}^{Hd}$, respectively, comprised of those supermanifolds which are locally modeled over superdomains in $\SDom_{MS}^{lko}$. 
\end{Defn}

\begin{Rem}
	To consider the DeWitt topology in the definition of supermanifolds is a little bit beside the point. It serves only to single out the correct model spaces. Indeed, by \thmref{Lem}{dersmooth}, all of the transition functions $U_{\alpha\beta}\to U_{\beta\alpha}$ are automatically homeomorphisms in the standard topology. 

	Tuynman \cite[Definitions 4.1]{tuynman} calls DeWitt--Tuynman supermanifolds without the Hausdorff axiom \emph{proto-$\mathcal A$-manifolds} (where $\mathcal A=\lambda^\infty$). 
\end{Rem}

\begin{Th}[dwtmsiso]
	Let $R$ be a non-discrete Hausdorff topological field of characteristic zero. The categories $\SMan_{dWT}$ and $\SMan_{MS}^{lko}$ are isomorphic. The isomorphism can be chosen to induce an isomorphism of $\SDom_{dWT}^{Hd}$ and $\SDom_{MS}^{lko,Hd}$. 
\end{Th}

\begin{proof}
	Let $\Phi:\SMan_{MS}^{lko}\to\SMan_{dWT}$ be defined on objects by $\Phi(\mathcal M,[\mathcal A])=(\varinjlim_N\mathcal M(\lambda^N),[\ger A])$, where $\ger A=(\varinjlim_N\vphi_{\alpha,\lambda^N})$ if $\mathcal A=(\vphi_\alpha)$, and on morphisms by $\Phi(f)=\varinjlim_Nf_{\lambda^N}$. It is fairly straightforward to check that $\Phi$ is well-defined and functorial.  
	
	Let $(M,[\ger A])\in\SMan_{dWT}$. Let $(\vphi_\alpha)\in[\ger A]$ be an atlas, and define \fa $\lambda\in\Lambda$ subspaces $\mathcal M(\lambda)=\bigcup_\alpha\vphi_\alpha(\mathcal U_\alpha(\lambda))$ of $M$ where $\mathcal U_\alpha\in\SDom_{MS}^{lko}$ is uniquely determined by the requirement that $\mathcal U_\alpha(\lambda^\infty)=U_\alpha=\dom\vphi_\alpha$. Then $\mathcal M(\lambda)$ is independent of the choice of atlas. 
	
	Moreover, $\mathcal M(\lambda)\approx\coprod_\alpha\mathcal U_\alpha(\lambda)/\sim$ as topological spaces, where $\mathcal U_{\alpha\beta}$ is determined by $\mathcal U_{\alpha,\beta}(\lambda^\infty)=\vphi_\alpha^{-1}(\vphi_\beta(U_\beta))$ and where $\sim$ is the equivalence relation on $\coprod_{\alpha,\beta}\mathcal U_{\alpha\beta}(\lambda)$ defined by the transition functions. Hence, there exist for any given even algebra morphism $\vrho:\lambda\to\lambda'$ continuous maps $\mathcal M(\vrho):\mathcal M(\lambda)\to\mathcal M(\lambda')$ determined by $\mathcal M(\vrho)\circ\vphi_\alpha|_{\mathcal U_\alpha(\lambda)}=\vphi_\alpha|_{\mathcal U_\alpha(\lambda')}\circ\mathcal U_\alpha(\vrho)$ \fa $\alpha$. Thus, we have defined a functor $\mathcal M\in\Top^\Lambda$. If we moreover define $\phi_{\alpha,\lambda}=\vphi_\alpha|_{\mathcal U_{\alpha}(\lambda)}$, then, by definition of $\mathcal M$, $\phi_\alpha=(\phi_{\alpha,\lambda})_\lambda$ is a natural transformation $\mathcal U_\alpha\Rightarrow\mathcal M$. 
	
	Furthermore, $\mathcal A=(\phi_\alpha:\mathcal U_\alpha\Rightarrow\mathcal M)$ is a covering in $\Top^\Lambda$. Since all of the maps $\vphi_{\beta\alpha}:U_{\alpha\beta}\to U_{\beta\alpha}$ are isomorphisms in $\SDom_{dWT}$, the fibred products $\mathcal U_\alpha\times_{\mathcal M}\mathcal U_\beta$ exist in $\SDom_{MS}^{lko}$, and it follows that $(\phi_\alpha)$ is an atlas. Moreover, equivalent (DeWitt--Tuynman-) atlases of $M$ lead by this construction to equivalent (Molotkov--Sachse-) atlases of $\mathcal M$. If, moreover, $\vphi:M\to N$ is a morphism of DeWitt--Tuynman supermanifolds, we obtain a natural transformation $\phi=(\phi_\lambda):\mathcal M\Rightarrow\mathcal N$ by setting $\phi_\lambda=\vphi|_{\mathcal M(\lambda)}$. By \thmref{Th}{dewitttuynman-ms-dom}, or directly by \thmref{Th}{tsmoothequiv}, this is well-defined. Now, let $\Psi(M,[\ger A])=(\mathcal M,[\mathcal A])$ and $\Psi(\vphi)=\phi$. This is easily seen to be a functor, and moreover, it is inverse to $\Phi$. 
	
	To show that $\SMan_{dWT}^{Hd}$ and $\SMan_{MS}^{lko.Hd}$ are isomorphic, one needs to see that $(M,[\ger A])$ is Hausdorff if and only if $\mathcal M(R)$ is. This follows in much the same way as in the proof of \thmref{Prop}{msglue}. 
\end{proof}

\bibliographystyle{alpha}%

\begin{thebibliography}{GGH10}

\bibitem[Ber87]{berezin}
F.A. Berezin.
\newblock {\em Introduction to superanalysis}, volume~9 of {\em Mathematical
  Physics and Applied Mathematics}.
\newblock D. Reidel Publishing Co., Dordrecht, 1987.

\bibitem[BGN04]{BGlN}
W.~Bertram, H.~Gl{\"o}ckner, and K.-H. Neeb.
\newblock Differential calculus over general base fields and rings.
\newblock {\em Expo. Math.}, 22(3):213--282, 2004.

\bibitem[Bou89]{bourbaki-top1}
N.~Bourbaki.
\newblock {\em General topology. {C}hapters 1--4}.
\newblock Elements of Mathematics. Springer-Verlag, Berlin, 1989.

\bibitem[DeW84]{dewitt}
B.S. DeWitt.
\newblock {\em Supermanifolds}.
\newblock Cambridge Monographs on Mathematical Physics. Cambridge University
  Press, Cambridge, 1984.

\bibitem[DM99]{deligne-morgan}
P.~Deligne and J.W. Morgan.
\newblock Notes on supersymmetry (following {J}oseph {B}ernstein).
\newblock In {\em Quantum fields and strings: a course for mathematicians,
  {V}ol. 1, 2 ({P}rinceton, {NJ}, 1996/1997)}, pages 41--97. Amer. Math. Soc.,
  Providence, RI, 1999.

\bibitem[FLV07]{fioresi-lledo-varadarajan}
R.~Fioresi, M.A. Lled{\'o}, and V.S. Varadarajan.
\newblock The {M}inkowski and conformal superspaces.
\newblock {\em J. Math. Phys.}, 48(11), 2007.

\bibitem[GGH10]{ggh-kac}
H.~Gl\"ockner, R.~Gramlich, and T.~Hartnick.
\newblock Final group topologies, {K}ac--{M}oody groups, and {P}ontryagin
  duality.
\newblock {\em Israel J. Math.}, 177:49--102, 2010.

\bibitem[Gir71]{giraud}
J.~Giraud.
\newblock {\em Cohomologie non ab\'elienne}, volume 179 of {\em Grundlehren der
  mathematischen Wissenschaften}.
\newblock Springer-Verlag, Berlin, 1971.

\bibitem[Gro60]{grothendieck-ega1}
A.~Grothendieck.
\newblock \'{E}l\'ements de g\'eom\'etrie alg\'ebrique. {I}. {L}e langage des
  sch\'emas.
\newblock {\em Inst. Hautes \'Etudes Sci. Publ. Math.}, 4:5--228, 1960.

\bibitem[Jas99]{jaskolski-smoothcoalg}
Z.~Jask{\'o}lski.
\newblock Smooth coalgebras.
\newblock {\em J. Geom. Phys.}, 29(1-2):87--150, 1999.

\bibitem[Kos77]{kostant-supergeom}
B.~Kostant.
\newblock Graded manifolds, graded {L}ie theory, and prequantization.
\newblock In {\em Differential geometrical methods in mathematical physics
  (Proc. Sympos., Univ. Bonn, Bonn, 1975)}, volume 570 of {\em Lecture Notes in
  Math.}, pages 177--306. Springer, Berlin, 1977.

\bibitem[Lei80]{leites}
D.~Leites.
\newblock Introduction to the theory of supermanifolds.
\newblock {\em Russian Math. Surveys}, 35(1):1--64, 1980.

\bibitem[Man88]{manin}
Y.I. Manin.
\newblock {\em Gauge field theory and complex geometry}, volume 289 of {\em
  Grundlehren der Mathematischen Wissenschaften}.
\newblock Springer-Verlag, Berlin, 1988.

\bibitem[Met03]{metzler}
D.S. Metzler.
\newblock Topological and smooth stacks.
\newblock arXiv:math/ 0306176v1, 2003.

\bibitem[Mol84]{molotkov}
V.~Molotkov.
\newblock Infinite-dimensional {$\mathbb Z_2^k$}-manifolds.
\newblock ICTP preprints IC/84/183, 1984.

\bibitem[Sa07]{sachse-diss}
C.~Sachse.
\newblock {\em Global analytic approach to super Teichm\"uller spaces}.
\newblock PhD thesis, Universit\"at Leipzig, 2007.

\bibitem[Sa08]{sachse-preprint}
C.~Sachse.
\newblock A categorical formulation of superalgebra and supergeometry.
\newblock arXiv:0802.4067, 2008.

\bibitem[SaW09]{sachse-wockel}
C.~Sachse and C.~Wockel.
\newblock The diffeomorphism supergroup of a finite-dimensional supermanifold.
\newblock arXiv:0904.2726v2 [math.DG], 2009.

\bibitem[Sch84]{schmitt-supergeom}
T.~Schmitt.
\newblock Super differential geometry.
\newblock Akademie der Wissenschaften der DDR, Institut f\"ur Mathematik,
  Report R-MATH-05/84, 1984.

\bibitem[Sch88]{schmitt-infdimsmf}
T.~Schmitt.
\newblock Infinite-dimensional supermanifolds. {I}.
\newblock Akademie der Wissenschaften der DDR, Institut f\"ur Mathematik,
  Report R-MATH-88-08, 1988.

\bibitem[Schw84]{schwarz}
A.S. Schwarz.
\newblock On the definition of superspace.
\newblock {\em Theoret. Math. Phys.}, 60(1):657--660, 1984.

\bibitem[Tuy97]{tuynman-ck}
G.M. Tuynman.
\newblock Functions of class $\mathcal C^k$ without derivatives.
\newblock Publ. Mat. 41(2):417--435, 1997.

\bibitem[Tuy04]{tuynman}
G.M. Tuynman.
\newblock {\em Supermanifolds and supergroups. Basic theory}, volume 570 of
  {\em Mathematics and its Applications}.
\newblock Kluwer Academic Publishers, Dordrecht, 2004.

\bibitem[Vis05]{vistoli}
A.~Vistoli.
\newblock Grothendieck topologies, fibered categories and descent theory.
\newblock In {\em Fundamental algebraic geometry}, volume 123 of {\em Math.
  Surveys Monogr.}, pages 1--104. Amer. Math. Soc., Providence, RI, 2005.

\bibitem[Vor84]{voronov-maps}
A.A. Voronov.
\newblock Mappings of supermanifolds.
\newblock {\em Theoret. Math. Phys.}, 60(1):660--664, 1984.

\bibitem[Yam98]{yamasaki}
A.~Yamasaki.
\newblock Inductive limit of general linear groups.
\newblock {\em J. Math. Kyoto Univ.}, 38(4):769--779, 1998.

\end{thebibliography}

\end{document}